\documentclass[12pt,twoside]{amsart}

\usepackage{amssymb,latexsym,bbm,graphicx,epsfig,epic,eepic,oldgerm,psfrag}
\usepackage{a4wide}

\usepackage{xypic}
\input xy

\xyoption{all}

\theoremstyle{plain}
\newtheorem{theorem}{Theorem}[section]
\newtheorem{lemma}[theorem]{Lemma}
\newtheorem{proposition}[theorem]{Proposition}
\newtheorem{corollary}[theorem]{Corollary}
\newtheorem*{remark*}{Remark}
\newtheorem*{remarks*}{Remarks}
\newtheorem{remark}[theorem]{Remark}

\newtheorem*{example*}{Example}
\newtheorem*{examples*}{Examples}

\newtheorem{question}[theorem]{Question}

\newcommand{\proofend}{\hspace*{\fill} $\Box$\\}
\newcommand{\diam}{\hspace*{\fill} $\Diamond$}

\newcommand{\Dcirc}{\overset%
{\raisebox{-.3ex}[0ex][-.3ex]{\mbox{$\scriptscriptstyle\circ$}}\mskip-5mu}D}

\def\1{\:\!}
\def\2{\;\!}
\def\s{\smallskip}
\def\m{\medskip}
\def\eps{\varepsilon}

\def\im{\operatorname {im}}

\def\Diffc0{\operatorname{Diff^c_0}}

\def\Sympc0{\operatorname{Symp^c_0}}

\def\rank{\operatorname{rank}}

\def\reg{\operatorname{reg}}

\def\top{\operatorname{top}}

\def\HF{\operatorname{HF}}
\def\CF{\operatorname{CF}}
\def\HM{\operatorname{HM}}
\def\CM{\operatorname{CM}}
\def\C{\operatorname{C}}
\def\H{\operatorname{H}}

\def\GL{\operatorname{GL}}
\def\SL{\operatorname{SL}}

\def\ga{\alpha}
\def\gb{\beta}
\def\gg{\gamma}
\def\gd{\delta}
\def\gve{\varepsilon}
\def\gf{\varphi}

\def\gl{\lambda}
\def\go{\omega}
\def\gs{\sigma}

\def\ca{{\mathcal A}}

\def\ce{{\mathcal E}}

\def\ch{{\mathcal H}}
\def\cj{{\mathcal J}}
\def\cl{{\mathcal L}}
\def\cm{{\mathcal M}}

\def\cp{{\mathcal P}}

\def\cs{{\mathcal S}}

\def\bJ{{\mathbf J}}

\def\FF{\mathbbm{F}}

\def\NN{\mathbbm{N}}
\def\PP{\mathbbm{P}}
\def\QQ{\mathbbm{Q}}
\def\RR{\mathbbm{R}}

\def\ZZ{\mathbbm{Z}}

\def\R{\operatorname{\mathbbm{R}}}
\def\RP{\operatorname{\mathbbm{R}P}}
\def\CP{\operatorname{\mathbbm{C}P}}

\def\pp{\partial}

\def\ra{\rightarrow}
\def\ha{\hookrightarrow}

\def\ni{\noindent}
\def\b{\bigskip}
\def\m{\medskip}

\def\proof{\noindent {\it Proof. \;}}
\begin{document}

\title{Positive topological entropy of Reeb flows on spherizations}

\author{Leonardo Macarini} 
\thanks{LM partially supported by CNPq}
\address{(L.~Macarini)
Universidade Federal do Rio de Janeiro, Instituto de Matem\'atica,
Cidade Universit\'aria, CEP 21941-909 - Rio de Janeiro - Brazil}
\email{leonardo@impa.br}
\author{Felix Schlenk} 
\thanks{FS partially supported by SNF grant 200021-125352/1.}
\address{(F.~Schlenk)
Institut de Math\'ematiques,
Universit\'e de Neuch\^atel,
Rue \'Emile Argand~11,
CP~158,
2009 Neuch\^atel,
Switzerland}
\email{schlenk@unine.ch}

\date{\today}

\begin{abstract}
Let $M$ be a closed manifold whose based loop space $\Omega (M)$ is ``complicated''.
Examples are rationally hyperbolic manifolds and manifolds whose fundamental group has exponential growth.
Consider a hypersurface $\Sigma$ in $T^*M$ which is fiberwise starshaped with respect to the origin.
Choose a function $H \colon T^*M \to \RR$ such that $\Sigma$ is
a regular energy surface of $H$,
and let $\gf^t$ be the restriction to $\Sigma$ of the
Hamiltonian flow of~$H$.

\m
\ni
{\bf Theorem~1.}
{\it The topological entropy of $\gf^t$ is positive.}

\m
\ni
This result has been known for fiberwise {\it convex}\,
$\Sigma$ by work of Dinaburg, Gromov, Paternain, and Paternain--Petean on geodesic flows.
We use the geometric idea and the Floer homological
technique from~\cite{FS3},
but in addition apply the sandwiching method.
Theorem~1 can be reformulated as follows.

\m
\ni
{\bf Theorem~1'.}
{\it The topological entropy of any Reeb flow on the spherization $SM$ of $T^*M$ is positive.}

\m
\ni
For $q \in M$ abbreviate $\Sigma_q = \Sigma \cap T_q^*M$.
The following corollary extends results of 
Morse and Gromov on the number of 
geodesics between two points.

\m
\ni
{\bf Corollary 1.}
{\it Given $q \in M$, for almost every $q' \in M$
the number of orbits of the flow $\gf^t$ from $\Sigma_q$ to $\Sigma_{q'}$ 
grows exponentially in time.
}

\m
\ni
In the lowest dimension,
Theorem~1 yields the existence of many {\it closed}\, orbits.

\m
\ni
{\bf Corollary 2.}
{\it Let $M$ be a closed surface different from $S^2$, $\RP^2$, the torus and the Klein bottle.
Then $\gf^t$ carries a horseshoe.
In particular, the number of geometrically distinct closed orbits
grows exponentially in time.
}
\end{abstract}

\maketitle


\section{Introduction and main results}

\ni
The topological entropy $h_{\top} (\gf)$ of a
diffeomorphism~$\gf$ of a compact manifold $P$ is a basic numerical
invariant measuring the orbit structure complexity of~$\gf$.
There are various ways of defining $h_{\text{top}} (\gf)$, 
see \cite{KH} and Section~\ref{s:vol.top} below.
In this paper, we show that for a certain class of diffeomorphisms $\gf$ on a certain class of manifolds~$P$, the number $h_{\text{top}} (\gf)$ is always positive.
We start with introducing the manifolds, first addressing the base manifolds $M$ and then the hypersurfaces $\Sigma \subset T^*M$, and then describe the diffeomorphisms studied.

\subsection{Energy hyperbolic manifolds}  \label{ss:energy}
The complexity of our flows on $\Sigma \subset T^*M$ comes from the complexity of the loop space of the base manifold~$M$.
Let $(M,g)$ be a $C^\infty$-smooth closed connected Riemannian manifold.
We assume throughout that $M$~is connected.
Fix $q_0 \in M$ and denote by $\Omega^1(M,q_0)$ the space of all paths $q \colon [0,1] \to M$ of Sobolev class $W^{1,2}$ such that $q(0)=q(1)=q_0$.
This space has a canonical Hilbert manifold structure, \cite{Kli}.
The energy functional $\ce = \ce_g \colon \Omega^1(M,q_0) \to \RR$ is defined as
\[
\ce (q) \,:=\, \frac 12 \int_0^1 | \dot q (t) |^2 \,dt
\]
where $| \dot q (t) |^2 = g_{q(t)} \left( \dot q (t) , \dot q (t) \right)$.
For $a > 0$ we consider the sublevel sets
\[
\ce^a(q_0) \,:=\, \left\{ q \in \Omega^1 ( M,q_0 ) \mid
\ce(q) \le a \right\} .
\]
Let $\PP$ be the set of non-negative integers which are prime or $0$, and write $\FF_p = \ZZ_p$ and $\FF_0 = \QQ$.
Throughout, $\H_*$ denotes singular homology.
Let
\[
\iota_k \colon \H_k \left( \ce^a(q_0); \FF_p \right) \,\to\,
                         \H_k \left( \Omega^1(M,q_0); \FF_p \right)
\]
be the homomorphism induced by the inclusion $\ce^a (q_0) \ha \Omega^1(M,q_0)$.
It is well-known that for each~$a$ the homology groups
$\H_k \left( \ce^a(q_0); \FF_p \right)$ vanish for all large enough~$k$, cf.~\cite{Benci}.
Therefore, the sums in the following definition are finite.
Following~\cite{FS3} we make the

\m
\ni
{\bf Definition} 
The Riemannian manifold $(M,g)$ is {\it energy hyperbolic}\,
if
$$
C(M,g) \,:=\,
\sup_{p \in \PP} \liminf_{m \to \infty} \frac 1m \log \sum_{k \ge 0}
\dim \iota_k \H_k \left( \ce^{\frac 12 m^2}(q_0) ; \FF_p \right) \,>\, 0.
$$

\ni
{\bf Remarks.}
Since $M$ is connected, $C(M,g)$ does not depend on $q_0$.
Since $M$ is closed, the property ``energy hyperbolic" does not depend on $g$.
We say that the closed manifold $M$ is {\it energy hyperbolic}\/ if $(M,g)$ is
energy hyperbolic for some and hence any Riemannian metric~$g$ on~$M$.
\diam

\b
\ni
{\bf Examples.} 
``Most'' closed manifolds are energy hyperbolic:

\m \ni 
{\bf 1.~Manifolds whose fundamental group has exponential growth.}
Assume that the fundamental group $\pi_1(M)$ has exponential growth.
This is, for instance, the case if~$M$ admits a Riemannian metric of negative sectional
curvature, \cite{Mi}.
By Proposition~1.8 of~\cite{FS3}, $M$ is energy hyperbolic.

\m \ni
{\bf 2.~Hyperbolic manifolds.}
Given a continuous map $f \colon (X,x) \to (Y, f(x))$ between path-connected pointed spaces,
let $\Omega(f) \colon \Omega (X) \to \Omega (Y)$ be the induced map between based loop spaces and let
$$
\H_* \left( \Omega (f); \FF_p \right) \colon
\H_* \left( \Omega (X);\FF_p\right) \to \H_* \left( \Omega (Y);\FF_p \right)
$$
be the map induced in homology.
Following~\cite{PP2} we say that a closed connected manifold~$M$
is {\it hyperbolic}\, if there exists a finite simply connected CW complex $K$ and a continuous map $f \colon K \to M$ such that for some $p \in \PP$ the sequence
$$
r_m (M,K,f; \FF_p) \,:=\,
\sum_{k=0}^m \rank \H_k \left( \Omega (f); \FF_p \right)
$$
grows exponentially in~$m$.
In particular, rationally hyperbolic manifolds 
(such as $\CP^2$ blown up in at least two points)
are hyperbolic, 
\cite{PP2}.
An example of a hyperbolic manifold that is not rationally hyperbolic is
$T^4 \# \overline{\CP}^2$.
By Proposition~1.11 of~\cite{FS3}, hyperbolic manifolds are energy hyperbolic.
We refer to~\cite{FS3} and the references therein for a thorough discussion of these two classes 
of energy-hyperbolic manifolds and for more examples.
\diam

\subsection{Fiberwise starshaped hypersurfaces in $T^*M$}

Let $\Sigma$ be a smooth connected hypersurface in $T^*M$. We say
that $\Sigma$ is {\it fiberwise starshaped}\, if for each point $q
\in M$ the set $\Sigma_q := \Sigma \cap T_q^*M$ is the smooth
boundary of a domain in $T_q^*M$ which is strictly starshaped with
respect to the origin $0_q \in T^*M$. This means that the radial
vector field $\sum_i p_i\;\!\pp p_i$ is transverse to each
$\Sigma_q$. We assume throughout that $\dim M \ge 2$. 
Then $T^*M \setminus \Sigma$ has two components, 
the bounded inner part $\Dcirc (\Sigma)$ containing the zero section, 
and the unbounded outer part $D^c(\Sigma) = T^*M \setminus D(\Sigma)$,
where $D(\Sigma)$ denotes the closure of $\Dcirc (\Sigma)$.

\subsection{Dynamics on fiberwise starshaped hypersurfaces}

Given a fiberwise starshaped hypersurface $\Sigma \subset T^*M$,
choose a smooth function $H \colon T^*M \to \RR$
such that $H^{-1}(1) = \Sigma$ and such that $1$ is a regular value of $H$.
Let
\begin{equation}  \label{e:go}
\go \,=\, \sum_{i=1}^d dp_i \wedge dq_i
\end{equation}
be the standard symplectic form on~$T^*M$.
The Hamiltonian vector field $X_H$ defined by
\begin{equation}  \label{e:XH}
\go(X_H, \cdot) \,=\, -dH(\cdot)
\end{equation}
defines a (local) flow $\gf_H^t$ on $T^*M$, called the Hamiltonian flow of $H$.
It restricts to a flow $\gf_H^t |_\Sigma$ on $\Sigma$.
Our sign convention in~\eqref{e:go} and~\eqref{e:XH} is such that
the flow $\gf_H^t$ of a geodesic Hamiltonian $H = \frac 12 |p|^2$ is the geodesic flow.

The orbits of $\gf_H^t|_\Sigma$ do not depend on~$H$ in the sense that 
for a different choice~$G$ of the Hamiltonian function,
the flow $\gf_G^t|_\Sigma$ is a time change of the flow $\gf_H^t|_\Sigma$,
i.e., $\gf_G^t|_\Sigma(x) = \gf_H^{\gs(t,x)}|_\Sigma(x)$ 
for a smooth positive function $\gs$ on $\RR \times \Sigma$.
For a flow without fixed points, 
any time change preserves vanishing of the topological entropy, see~\cite[p.~113]{KH}.
Therefore $h_{\top} (\gf_H^t |_\Sigma) >0$ iff $h_{\top} (\gf_G^t |_\Sigma) >0$.
Similarly, since $\Sigma$ is compact, given $q,q' \in M$ the number of $\gf_H^t$-orbits
from $\Sigma_q$ to $\Sigma_{q'}$ grows exponentially in time iff 
this is so for the number of $\gf_G^t$-orbits.

\subsection{Main result}

$\phantom{...}$

\m
\ni
{\bf Theorem 1.}
{\it Consider an energy hyperbolic manifold $M$ and a fiberwise starshaped hypersurface $\Sigma \subset T^*M$. Then $h_{\top} (\gf_H^t |_\Sigma) >0$.}

\b \ni
{\bf Discussion.}
{\bf 1.} 
The unit cosphere bundle
$$
S_1M (g) \,:=\, \left\{ (q,p) \in T^*M \mid \left| p \right| =1 \right\}
$$
associated to a Riemannian metric~$g$ on~$M$
is an example of a fiberwise starshaped hypersurface.
In this special case, Theorem~1 has been proved by Dinaburg~\cite{Dina}, Gromov~\cite{Gromov78},
Paternain~\cite{Pa1,Pa}, and Paternain--Petean~\cite{PP2} in their study of geodesic flows.
More generally, the sets $\Sigma_q$ are all convex if and only if $\Sigma$ is
the level set of a Finsler metric.
In this case, Theorem~1 also follows from Paternain's work,~\cite{Pa}.

\m
{\bf 2.}
The assumption that $M$ is energy hyperbolic obviously cannot be omitted:  
Geodesic flows over round spheres or flat tori have vanishing topological entropy.
The assumption that $\Sigma$ is fiberwise starshaped {\it with respect to the origin}\/
cannot be omitted either: 
There are exact magnetic flows on compact quotients of {\bf Sol} 
whose topological entropy drops to zero when the energy levels cease 
to inclose the zero section, see Section~\ref{s:example}.

\m
{\bf 3.}
Recall that a hypersurface $\Sigma \subset T^*M$ is said to be
{\it of restricted contact type}\/
if there exists a vector field $X$ on $T^*M$ such that $\cl_X \go = d \iota_X \go = \go$
and such that $X$ is everywhere transverse to $\Sigma$, pointing outwards.
Equivalently, there exists a $1$-form $\ga$ on $T^*M$ such that $d
\ga = \go$ and such that $\ga \wedge (d \ga)^{d-1}$ is a volume form
on $\Sigma$ orienting $\Sigma$ as the boundary of~$D(\Sigma)$. 
The correspondence is given by $\ga = \iota_X \go$. 
Our assumption that $\Sigma$ is fiberwise
starshaped translates to the assumption that $\Sigma$ is of
restricted contact type {\it with respect to}\/ the Liouville vector
field $Y = \sum_{i=1}^d p_i \2 \pp p_i$, 
or, equivalently, the Liouville form
$$
\gl \,=\, \sum_{i=1}^d p_i \2 dq_i .
$$
defines a contact form on $\Sigma$.

\begin{question}
{\rm
Is Theorem 1 true for any hypersurface $\Sigma \subset T^*M$ 
of restricted contact type that encloses the $0$-section?}
\end{question}

\subsection{Reformulation of Theorem~1}

We are going to reformulate Theorem~1,
which is formulated in terms of Hamiltonian dynamics,
in a more invariant way by means of contact geometry.
Let~$M$ be a closed connected manifold.
Fix a fiberwise starshaped hypersurface $\Sigma \subset T^*M$.
The hyperplane field $\xi_{\Sigma} = \ker \left( \gl |_\Sigma \right) \subset T\Sigma$
is a contact structure on $\Sigma$.
If~$\Sigma'$ is another fiberwise starshaped hypersurface, then the differential 
of the diffeomorphism
$$
\Psi_{\Sigma \Sigma'} \colon \Sigma \to \Sigma', \quad 
(q,p) \mapsto \bigl( q, \psi(q,p) \2 p \bigr) ,
$$
obtained by fiberwise radial projection,
maps $\xi_\Sigma$ to $\xi_{\Sigma'}$ 
(because $\Psi_{\Sigma \Sigma'}^* \left( \gl |_{\Sigma'} \right) = \psi \1 \gl |_\Sigma$), 
and is hence a contactomorphism 
$\left( \Sigma, \xi_{\Sigma} \right) \to \left( \Sigma', \xi_{\Sigma'} \right)$.
These contact manifolds can thus be identified, and are called 
the {\it spherization} $(SM, \xi)$ of the cotangent bundle $(T^*M, \go)$.
\footnote{The authors of \cite{EKP} call it
{\it space of oriented contact elements} and write $\PP_+T^*M$.
We have chosen the above wording since the unit cosphere bundle
$(S_1M(g), \ker \gl)$ with respect to a Riemannian metric $g$ is a representative.}

Fix a representative $\left( \Sigma, \xi_{\Sigma} \right)$.
For every positive smooth function $f \colon \Sigma \to \RR$ we still have 
$\xi_{\Sigma} = \ker \left( f \1 \gl |_{\Sigma} \right)$.
The {\it Reeb vector field}\/ $R_f$ on $T\Sigma$ is defined as the unique vector field
such that
$$
d(f\gl) (R_f, \cdot ) \equiv 0, \qquad f \1 \gl (R_f) \equiv 1 .
$$
Its flow is called the {\it Reeb flow}\/ of $R_f$.
For $f \equiv 1$ the Reeb flow of $R_f$ is a time change of the flow
$\gf_H^t |_{\Sigma}$ of any Hamiltonian function $H \colon T^*M \to \RR$ with
$H^{-1}(1) = \Sigma$ and such that $1$ is a regular value.
For different functions $f$ the Reeb flows on $\Sigma$ can be completely different,
see e.g.~\cite{HMS}.

Given another hypersurface $\Sigma'$, let $\psi \colon \Sigma \to \RR$ be the positive function 
such that
$$
\Sigma' \,=\, \left\{ \bigl( q, \psi (q,p)\2 p \bigr) \mid (q,p) \in \Sigma \right\} .
$$
Then $d \2 \Psi_{\Sigma \Sigma'} (R_\gl) = R_{\psi^{-1} \gl}$, that is,
$\Psi_{\Sigma \Sigma'}$ conjugates the Reeb flows on 
$\bigl( \Sigma, \gl \bigr)$ and $\bigl( \Sigma', \psi^{-1}\gl \bigr)$.
Summarizing, we have that the set of Reeb flows on $(SM,\xi)$ is in bijection with Hamiltonian flows
on fiberwise starshaped hypersurfaces in $T^*M$, up to time changes.
Recall that for a flow without fixed points, 
any time change preserves vanishing of the topological entropy. 
Theorem~1 is therefore equivalent to 

\m \ni 
{\bf Theorem 1'.}
{\it 
The topological entropy of any Reeb flow on the spherization $SM$ of $T^*M$ is positive.
}

\subsection{Uniform exponential growth of the number of Reeb chords}

Let $M$, $\Sigma$ and $H$ be as in Theorem~1.
For $q,q' \in M$ let
$\nu_T (q,q',H)$
be the number of flow lines of $\gf_H^t|_\Sigma$
starting from $\Sigma_q$ at $t=0$ and arriving on $\Sigma_{q'}$
before time~$T$.
From our proof of Theorem~1 we shall immediately obtain

\m
\ni
{\bf Corollary 1.}
{\it 
There exists a constant $h>0$ 
depending only on $(M,g)$, $\Sigma$ and~$H$
such that for each
$q \in M$,
$$
\liminf_{n \to \infty} \frac 1n \log \nu_n(q,q',H) \,\ge\, h
$$
for almost every $q' \in M$.
Moreover, 
for every $q'$ there exists an orbit of $\gf_H^t$ from $\Sigma_q$ 
to~$\Sigma_{q'}$.
}

\m
For geodesic flows such lower bounds were obtained by
Morse and Gromov~\cite{Gromov78}.
%
%
Denote the fiber of $SM$ over $q$ by $S_qM$.
In terms of Reeb flows, Corollary~1 reads

\m
\ni
{\bf Corollary 1'.}
{\it 
Fix a Reeb flow $\gf^t$ on $(SM,\xi)$.
There exists a constant $h>0$ depending only on $(M,g)$ and $\gf^t$
such that for each
$q \in M$,
$$
\liminf_{n \to \infty} \frac 1n \log \nu_n(q,q',\gf^t) \,\ge\, h
$$
for almost every $q' \in M$.
Moreover,
for every $q'$ there exists a Reeb chord from $S_qM$ to~$S_{q'}M$.
}

\m
Note that $S_qM$ and $S_{q'}M$ are Legendrian submanifolds of $\bigl( SM,\xi \bigr)$.
Corollary~1' is thus a special case of the Arnold chord conjecture, with multiplicities.

\m
Let $(M,g)$ be a closed Riemannian manifold, endowed with its 
Riemannian measure~$\mu_g$,
and let $\gf_g$ be the geodesic flow on $S_1M(g)$.
It has been shown by 
Ma\~n\'e \cite{Mane} and Paternain--Paternain~\cite{PaPa} that 
\begin{equation} \label{e:MPP}
h_{\top}(\gf_g) \,=\, \lim_{T \to \infty} \frac 1T \log 
\int_{M \times M} \nu_T(q,q') \2 d\mu_g(q) \1 d\mu_g(q') .
\end{equation}

\begin{question}
{\rm
Is the identity~\eqref{e:MPP} true 
for all Reeb flows on~$SM$?
}
\end{question}

\ni
For a partial answer in the positive in the case of certain exact magnetic flows see~\cite{Niche}.

\m
In dimension~2, Theorem~1 yields the existence of many {\it closed}\, orbits.
Indeed, according to \cite{Katok.IHES} and \cite[Theorem~4.1]{Katok.ETDS},
a smooth flow on a closed $3$-manifold with positive topological entropy
has exponential growth of closed orbits. 

\m \ni
{\bf Corollary 2.}
{\it Let $M$ be a closed surface different from $S^2$, $\RP^2$, the torus and the Klein bottle.
Then $\gf_H^t$ carries a horseshoe.
In particular, the number of geometrically distinct closed orbits
grows exponentially in time.
}

\m
This is a special case of the Weinstein conjecture, with multiplicities.
In terms of Reeb flows, Corollary~2 means that for such surfaces~$M$, 
{\it for any Reeb flow on $SM$ the number of geometrically distinct closed Reeb chords 
grows exponentially in time.}
For a generalization of this result to generic Reeb flows on spherizations over higher dimensional closed manifolds~$M$
we refer to \cite{Heistercamp}.
For a result on exponential growth rate of the number of closed orbits for certain Reeb flows
on a large class of closed contact $3$-manifolds, see~\cite{CoH}.

\begin{remark}
{\rm
In this paper we have exploited the {\it exponential}\/ growth
of the homology of the based loop space to obtain lower bounds for the entropy 
and the number of Reeb chords.
For base manifolds~$M$ whose based loop space grows only polynomially,
the same method yields lower bounds for the {\it slow entropy}, and {\it polynomial}\/
lower bounds for the number of Reeb chords, see~\cite{Wullschleger}.
In particular, for any two points $q,q'$ in a closed manifold~$M$ and any Reeb flow on~$SM$
there exists a Reeb chord from $\Sigma_q$ to $\Sigma_{q'}$.
\diam 
}
\end{remark}


\m
The paper is organized as follows.
In the next section we give the idea of the proof.
In Section~\ref{s:Ham.spec} 
we define the relevant Hamiltonians and compute their action spectra.
In Section~\ref{s:Floer} we recall Lagrangian Floer homology and compute it for 
our two Lagrangians in question.
In Sections~\ref{s:th1} and \ref{s:c1} we derive Theorem~1 and Corollary~1.
Finally, in Section~\ref{s:example} we discuss the exact magnetic flows mentioned in 
Discussion~1.2.

\m
\ni
{\bf Acknowledgments.}
We cordially thank 
Albert Fathi, Emmanuel Giroux, 
Evgenij Gutkin, Anatole Katok, 
Gabriel Paternain, Leonid Polterovich and Otto van Koert
for valuable discussions. 
Special thanks to Gabriel Paternain for explaining to us the example in Section~\ref{s:example}. 
Most of this paper was written during the second author's visit at IMPA in February and August 2005 and in January 2008, and the first author's visit at ULB in September 2007.
We thank both institutions for their warm hospitality.

\section{Idea of the proof}

Let $M$ be an energy-hyperbolic manifold, 
and let $\Sigma \subset T^*M$ be a fiberwise starshaped hypersurface.
Recall that $D(\Sigma)$ denotes the closure of the bounded part of
$T^*M \setminus \Sigma$. 
For $q \in M$ let
$$
D_q(\Sigma) \,=\, D(\Sigma) \cap T^*_qM
$$
be the closed starshaped disc over $q$ with boundary $\Sigma_q = \Sigma \cap T^*_qM$.
Choose a smooth Hamiltonian function $K \colon T^*M \to \RR$
that is fiberwise homogenous of degree~2, 
up to a smoothening near the zero section~$M$.
Recall that it suffices to show that $h_{\top} (\gf_K |_{\Sigma}) >0$.
Since $K$~is fiberwise homogenous of degree~2, its Hamiltonian flows
on the levels $s \1 \Sigma$ agree up to constant time changes.
Therefore, $h_{\top} (\gf_K |_{\Sigma}) >0$ is equivalent to 
$h_{\top} (\gf_K |_{D (\Sigma)}) >0$.
By a theorem of Yomdin, the latter inequality will follow 
if we can find a point $q \in M$ such that 
the volume of the discs $\gf_K^n \bigl( D_q(\Sigma) \bigr)$ 
grows exponentially with~$n$.

We shall, in fact, establish this for each point $q \in M$,
by taking up an idea from~\cite{FS:GAFA}: 
Fix $q \in M$.
By sandwiching the set $\Sigma$ between two sphere bundles
(which are the levels of a geodesic Hamiltonian),
and by using the Lagrangian Floer homology of two fibers in~$T^*M$
and its isomorphism to the homology of the based loop space of~$M$
(whose dimensions grow exponentially by assumption),
we shall show that for almost every other point $q' \in M$
the number of intersections
$$
\gf_K^n \bigl( D_q(\Sigma) \bigr) \cap D_{q'}(\Sigma)
$$
grows exponentially with $m$. 
This means that the discs $\gf_K^n \bigl( D_q(\Sigma) \bigr)$
``wrap'' exponentially often in~$n$ around the base~$M$,
and hence their volume grows exponentially.

\begin{remark}
{\rm
Our proof of Theorem~1 can be considerably shortened by using
symplectic homology with its wrapped version for Lagrangian boundary
conditions and a remark of Seidel in~\cite[Section 4a]{Seidel.biased},
see~\cite{MMP}.
Furthermore, it is conceivable that Theorem~1 (or at least Corollary~1) 
can also be obtained by using Contact homology,
for spherizations, relative to two Legendrian submanifolds $\Sigma_q$ and $\Sigma_{q'}$,
or by a relative version of Rabinowitz--Floer homology,
and by establishing an isomorphism between either of these homologies with 
the homology of the based loop space.
This would have the advantage of working directly at the level~$\Sigma$, 
instead of the sublevel set~$D(\Sigma)$.
These homologies are, however, not yet fully established. 
We thus preferred to work with Lagrangian Floer homology and to pass through the sublevel set~$D(\Sigma)$.
} 
\end{remark}

\section{Relevant Hamiltonians and their action spectra}  \label{s:Ham.spec}

\subsection{The path space $\Omega^1 \left( T^*M,q_0,q_1 \right)$}
For $q_0, q_1 \in M$ let $\Omega^1 \left(T^*M,q_0,q_1\right)$ be the
space of all paths $x \colon [0,1] \to T^*M$ of Sobolev class
$W^{1,2}$ such that $x(0) \in T_{q_0}^*M$ and $x(1) \in T_{q_1}^*M$.
This space has a canonical Hilbert manifold structure, \cite{Kli}.
Consider a proper Hamiltonian function $H \colon T^*M \to \RR$; then
its Hamiltonian flow $\gf_H^t$ exists for all times. The action
functional of classical mechanics $\ca_H \colon \Omega^1 \left(
T^*M,q_0,q_1 \right) \to \RR$ associated with $H$ is defined as
\begin{equation}  \label{def:af}
\ca_H(x) \,=\, \int_0^1 \bigl( \gl \left(\dot{x}(t) \right) -
H(x(t)) \bigr) dt,
\end{equation}
where $\gl = \sum_{j=1}^d p_j \,dq_j$ is the canonical $1$-form on
$T^*M$. This functional is $C^\infty$-smooth, and its critical
points are precisely the elements of the space $\cp \left( H,q_0,q_1 \right)$ 
of $C^\infty$-smooth paths $x \colon [0,1] \to T^*M$
solving
$$
\dot{x}(t) = X_H \bigl( x(t) \bigr), \,\,t \in [0,1], \quad \; x(j)
\in T_{q_j}^*M, \,\,j=0,1.
$$
Notice that the elements of $\cp \left( H,q_0,q_1 \right)$
correspond to the intersection points of $\gf_H \big( T^*_{q_0}M \bigr)$ 
and $T^*_{q_1}M$ via the evaluation map $x \mapsto x(1)$.

\subsection{The Hamiltonians $G_- \le K \le G_+$}

\begin{figure}[h]
 \begin{center}
  \psfrag{e}{$\gve$}
  \psfrag{e2}{$\gve^2$}
  \psfrag{r}{$r$}
  \psfrag{f}{$f(r)$}
  \leavevmode\epsfbox{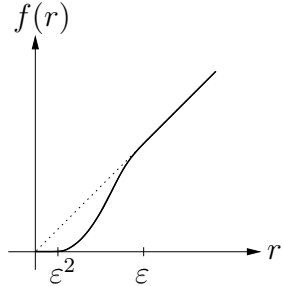}
 \end{center}
 \caption{The ``cut off'' function~$f$.}
 \label{figure.cutoff}
\end{figure}

Since $\Sigma$ is fiberwise starshaped, we can define a function $F
\colon T^*M \to \RR$ by the two requirements
$$
F|_\Sigma \equiv 1, \qquad F(q,sp) = s^2 F(q,p) \quad \text{for all
$s \ge 0$ and $(q,p) \in T^*M$}.
$$
This function is fiberwise homogenous of degree~2, of class $C^1$,
and smooth off the zero-section. To smoothen $F$, choose $\gve \in
(0, \frac 14)$. We shall appropriately fix $\eps$ later on. Choose a
smooth function $f \colon \RR \to \RR$ such that
\begin{equation} \label{d:f}
\left\{
\begin{array} {ll}
f(r) = 0 & \text{ if }  r \le \gve^2, \\ [0.3em]
f(r) = r & \text{ if }  r \ge \gve, \\[0.3em]
f'(r) >0 & \text{ if }  r > \gve^2, \\[0.3em]
0 \le f'(r) \le 2 & \text{ for all } r,
\end{array}
\right.
\end{equation}
see Figure~\ref{figure.cutoff}.

%
%

\ni 
Then $f \circ F$ is smooth. 
The Lagrangian Floer homology of~$f \circ F$ (with respect to two fibers of $T^*M$) is, 
in general, not defined, since the moduli spaces of Floer trajectories may not be
compact. We shall therefore alter $f \circ F$ near infinity to a
Riemannian Hamiltonian. 
We denote canonical coordinates on $T^*M$ by $(q,p)$. 
Fix a Riemannian metric~$g$ on~$M$, 
and let $| \:\! \cdot \:\! |$ be the norm on the fibers of $T^*M$ induced by~$g$.
Define $G(q,p) = \frac 12 \left| p \right|^2$. 
Recall that $M$ is compact and that $\Sigma$ is fiberwise starshaped. 
After multiplying~$g$ by a positive constant, we can assume that
$G \le F$.
Choose $\gs >0$ such that $\gs G \ge F$.

%
For $r>0$ we abbreviate
$$
D(r) \,=\, \left\{ (q,p) \in T^*M \mid |p| \le r \right\} .
$$
Note that $D(2) = \{ G \le 2 \}$.
Choose a smooth function $\tau \colon \RR \to \RR$ such that
\begin{eqnarray*}
 \begin{array}{lll}
   \tau(r) = 0 & \text{ if } & r \le 2, \\ [0.3em]
   \tau(r) = 1 & \text{ if } & r \ge 4, \\[0.3em]
   \tau'(r) \ge 0 & \text{ for all } r \in \RR.
 \end{array}
\end{eqnarray*}
Define
\begin{eqnarray}
G_+    &=& \gs G \label{def:G+}\\
K(q,p) &=& \bigl( 1-\tau(|p|) \bigr) (f \circ F)(q,p) + \tau (|p|)
G_+(q,p),  \label{def:K} \\
G_-(q,p) &=& \bigl( 1-\tau(|p|) \bigr) (f \circ G)(q,p) + \tau
(|p|) G_+(q,p),  \label{def:G-}
\end{eqnarray}
Then
\begin{equation}  \label{e:KGG}
G_- \le K \le G_+
\end{equation}
and
\begin{equation}  \label{e:Dr}
K = f \circ F \;\text{ and }\; G_- = f \circ G \quad \text{on }
D(2) .
\end{equation}
Since $D(\Sigma) = \{ F \le 1 \} \subset \{ F \le 2 \} \subset \{ G_-
\le 2 \}$, we in particular have
\begin{equation}  \label{e:fSigma}
K = f \circ F  \quad \text{on } D(\Sigma) .
\end{equation}
Moreover,
\begin{equation}
G_- = K = G_+ \quad \text{ outside } D(4).
\end{equation}
\begin{figure}[h]
 \begin{center}
  \psfrag{G=}{$G=\frac 12 |p|^2$}
  \psfrag{G}{$G$}
  \psfrag{G-}{$G_-$}
  \psfrag{G+}{$G_+$}
  \psfrag{K}{$K$}
  \psfrag{2}{$2$}
  \psfrag{8}{$8$}
  \leavevmode\epsfbox{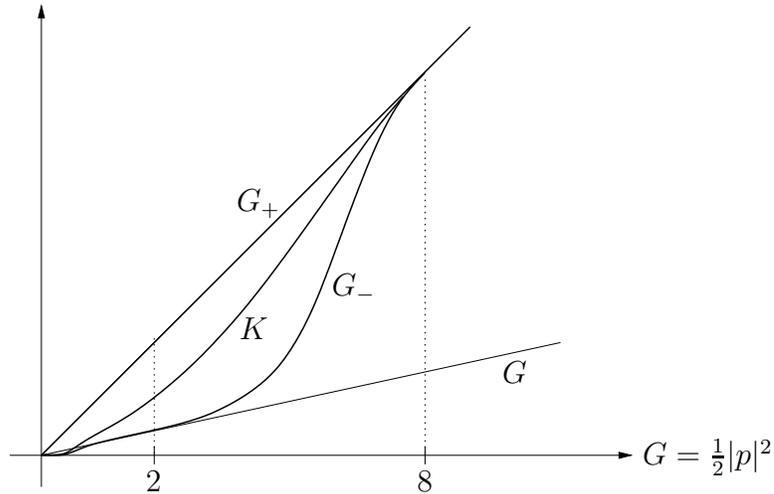}
 \end{center}
 \caption{The functions~$G_- \le K \le G_+$, schematically.}
 \label{figure.GKG}
\end{figure}
%
%

\subsection{Action spectra of $nG_-$, $nK$, $nG_+$}
We next investigate the action spectra of the functions $nG_-$, $nK$, $nG_+$ 
for $n \in \NN$. 
The action spectrum $\cs (H,q_0,q_1)$ of a proper Hamiltonian
$H \colon T^*M \to \RR$ is the set of critical values of $\ca_H
\colon \Omega^1(T^*M,q_0,q_1) \to \RR$,
\[
\cs (H,q_0,q_1) \,:=\, \left\{ \ca_H(x) \mid x \in \cp(H,q_0,q_1)
\right\} .
\]

We first look at the action spectra of Hamiltonians homogenous of
degree~2.

\begin{lemma}  \label{l:action}
Let $H \colon T^*M \to \RR$ be fiberwise homogenous of degree~2, let
$h \colon \RR \to \RR$ be a smooth function, and let $c>0$.

\s
\begin{itemize}
\item[(i)] 
For $\gg \in \cp (h \circ H,q_0,q_1)$ we have
$$
\ca_{h \circ H} (\gg) \,=\, 2\:\!h'\bigl(H(\gg)\bigr)
\:\!H(\gg)-h\bigr(H(\gg)\bigr).
$$

\s
\item[(ii)] 
$\cs (cH,q_0,q_1) = \frac 1c \cs (H,q_0,q_1)$.
\end{itemize}
\end{lemma}

\proof 
(i)
Let $Y = \sum_i p_i \;\! \pp p_i$. For each $t \in [0,1]$
we have $\dot \gg (t) = X_{h \circ H} \left(\gg(t)\right)$, and hence
$$
\gl \left( \dot \gg (t) \right) \,=\, \go \bigl( Y, X_{h \circ H}
(\gg(t))\bigr) \,=\, d(h \circ H) \left( \gg (t) \right) (Y) \,=\,
h' \bigl( H(\gg (t)) \bigr) dH ( \gg (t) ) (Y).
$$
Since $H$ is fiberwise homogenous of degree~2, Euler's identity
yields
$$
dH (\gg (t)) (Y) \,=\, 2\:\!H (\gg (t))
$$
whence
$$
\gl (\dot \gg (t)) \,=\, 2h' (H(\gg(t)) H(\gg(t)).
$$
Therefore,
$$
\ca_{h \circ H}(\gg) \,=\, \int_0^1 \bigl( \gl (\dot \gg(t))-
h(H(\gg(t))) \bigr) dt \,=\, 2 h'(H(\gg)) H(\gg) - h(H(\gg)),
$$
as claimed.

\s 
(ii) By definition~\eqref{e:XH} we have $X_{cH}(q,p) = c X_H
(q,p)$, and since $H$ is homogenous of degree~2, we have $X_H (q,
\frac 1c p) = \frac 1c X_H (q,p)$ (where in the latter equality we
have canonically identified $T_{(q,\frac 1c p)}T^*M =
T_{(q,p)}T^*M)$). Therefore,
\[
X_{cH}(q,\tfrac 1c p) \,=\, c X_H (q, \tfrac 1c p) \,=\, X_H(q,p) .
\]
To $\gg (t) = (q(t),p(t))$ in $\cp (H,q_0,q_1)$ therefore
corresponds $\gg_c (t) = \bigl( q(t), \frac 1c p(t) \bigr)$ in $\cp
(cH,q_0, q_1)$. Assertion~(i) now yields
\[
\ca_{cH} (\gg_c) \,=\, cH (\gg_c) \,=\, c H (q(0), \tfrac 1c p(0))
\,=\, \tfrac 1c H( q(0),p(0) ) \,=\, \tfrac 1c \ca_H (\gg),
\]
and so $\cs (cH,q_0,q_1) = \frac 1c \cs (H,q_0,q_1)$. \proofend

We next look at the action spectrum of $nK$.

\begin{proposition}  \label{p:action}
Let $\gg \in \cp (nK,q_0,q_1)$. 
If $\gg \subset \Dcirc(\Sigma)$, then
$\ca_{nK} (\gg) < n$. 
If $\gg \subset D^c (\Sigma)$, 
then $\ca_{nK} (\gg) > n$.
\end{proposition}

\proof 
Assume first that $\gg \subset \Dcirc(\Sigma)$. 
Then $F(\gg) < 1$ and $\gg \in \cp (n(f \circ F),q_0,q_1)$ by~\eqref{e:fSigma}.
Hence Lemma~\ref{l:action}\:(i) and our choice of~$f$ yield
$$
\ca_{nK}(\gg) \,=\, \ca_{n(f \circ F)}(\gg) \,=\, 2nf'(F(\gg))
F(\gg) - nf(F(\gg)) \,<\, n.
$$
Assume now that $\gg \subset D^c(\Sigma)$. 
Then $F(\gg) > 1$ and $n(f \circ F)(\gg) = nF(\gg)$. 
Let again $Y = \sum_i p_i \;\! \pp p_i$. 
For $t \in [0,1]$ we compute at $\gg (t) = (q,p)$, using~\eqref{def:K},
\begin{eqnarray*}
d(nK) (\gg (t)) (Y)
&=& - \tau'(|p|) nF(q,p) |p| + \left( 1-\tau(|p|) \right) d(nF)(q,p)(Y) \\
& & + \tau'(|p|) nG_+(q,p) |p| + \phantom{1-\2} \, \tau(|p|)\,
d(nG_+)(q,p)(Y).
\end{eqnarray*}
Since $\tau'(|p|) \ge 0$ and $nG_+ (q,p) \ge nF(q,p)$, the sum of
the first and third summand is $\ge 0$. Together with Euler's
identity applied to the functions $nF$ and $nG_+$ fiberwise
homogenous of degree~2 we obtain
\begin{eqnarray*}
d(nK) (\gg (t)) (Y) - nK(\gg (t)) &\ge&
\bigl( 1-\tau(|p|) \bigr) \bigl( d(nF)(\gg(t)) (Y) - nF(\gg(t)) \bigr) \\
&&
\phantom{1\;\,}+ \tau (|p|) \bigl( d(nG_+)(\gg(t)) (Y) - nG_+(\gg(t)) \bigr) \\
&=&
\left( 1-\tau(|p|) \right) nF(\gg(t)) + \tau(|p|) nG_+(\gg(t)) \\
&\ge& nF(\gg(t)) \\
&>& n.
\end{eqnarray*}
Hence $\ca_{nK} (\gg) = \int_0^1 \bigl( d(nK) (\gg (t)) (Y) - nK
(\gg (t)) \bigr) dt > n$, as claimed. \proofend

\subsubsection{The Non-crossing lemma}

For $n \in \NN$ consider the space of Hamiltonian functions
\[
\ch_4(nG_+) \,=\, \left\{ H \colon T^*M \to \RR \mid H = nG_+
\text{ on } T^*M \setminus D(4) \right\}.
\]
Note that $nG_-$, $nK$, $nG_+$ belong to $\ch_4(nG_+)$.
For $a \in \RR$ define
$$
\ch_4^a(nG_+) \,=\, \left\{ H \in \ch_4 (nG_+) \mid a \notin
\cs(H,q_0,q_1) \right\} .
$$
Fix a smooth function $\gb \colon \RR \to [0,1]$ such that
\begin{equation}  \label{e:beta}
\gb (s) =0 \,\text{ for } s \le 0, \quad \gb(s) =1 \,\text{ for } s \ge
1, \quad \gb'(s) \ge 0 \,\text{ for all } s \in \RR .
\end{equation}
For $s \in [0,1]$ define the functions
\begin{equation}  \label{def:Gs}
G_s \,=\, \bigl(1-\gb(s)\bigr) G_- + \gb (s) G_+ .
\end{equation}
Then $nG_s \in \ch_4(nG_+)$ for each $n \in \NN$ and $s \in
[0,1]$.

Recall that $\gs \ge 1$.
Choose $a \in (n,n+1)$ and define the
function $a(s) \colon [0,1] \to \RR$ by
\begin{equation}  \label{d:fas}
a(s) \,=\, \frac{a}{1+\gb(s)(\gs-1)} .
\end{equation}
Note that $a(s)$ is monotone decreasing with minimum $a(1) = a/\gs$.
We assume from now on that $\gve \in (0,\frac 14)$ entering the
definition of the smoothing function $f$ also meets
\begin{equation}  \label{ie:gve}
\gve^2 < \frac{1}{2\gs}.
\end{equation}
\begin{lemma}  \label{l:NC}
If $a \notin \cs (nG,  q_0,q_1)$, then 
$a(s) \notin \cs (nG_s,q_0,q_1)$
for $s \in [0,1]$.
\end{lemma}

\proof 
Take $\gg \in \cp(n G_s,q_0,q_1)$. In view of~\eqref{def:G-}
and \eqref{def:Gs}, the orbit $\gg$ lies on a level set of~$G$.

\s Assume first that $\gg \in \cp(nG_s)$ lies outside $D(2)$.
Then $f \circ G (\gg) = G(\gg)$, so that
\[
nG_s \,=\, n\bigl(1-\gb(s)\bigr) \Bigl( \bigl(1-\tau(|p|)\bigr)
G+\tau(|p|) G_+ \Bigr) + n \gb(s)G_+ .
\]
Computing as in the proof of Proposition~\ref{p:action} 
(with $F$ replaced by~$G$) we find
\begin{eqnarray*}
d(nG_s)(\gg(t)) (Y) - nG_s(\gg(t)) &\ge& n \bigl(1-\gb(s)\bigr)
\Bigl( \bigl( 1-\tau(|p|)\bigr) G(\gg(t)) + \tau(|p|) G_+(\gg(t))
\Bigr) \\
&& \phantom{n \bigl(1-\gb(s)\bigr) \bigl( G(\gg(t))\bigr)}+ n \gb(s)G_+(\gg(t)) \\
&\ge& n \bigl(1-\gb(s)\bigr) \bigl( G(\gg(t))\bigr)
+n\gb(s)G(\gg(t)) \\
&=& n G(\gg(t)) \\
&\ge& 2n .
\end{eqnarray*}
Hence $\ca_{nG_s}(\gg) \ge 2n \ge n+1 > a \ge a(s)$.

\s 
Assume next that $\gg$ lies in $D(\gve)$. 
Then $\tau(|p(\gg(t))|) = 0$ for all $t \in [0,1]$. 
Hence $G_-(\gg) = f \circ G(\gg)$
and
\begin{eqnarray*}
G_s(\gg) &=& \bigl(1-\gb(s) \bigr) \bigl( f \circ G \bigr)(\gg)
+\gb(s) \gs G(\gg) . 
\end{eqnarray*}
Applying Lemma~\ref{l:action}\:(i) with $h = \bigl( 1-\gb(s) \bigr) f + \gb(s) \gs$
we therefore find
\begin{eqnarray*}
\ca_{nG_s}(\gg) &=& 2n \Bigl( \bigl(1-\gb(s) \bigr) f'\bigl( 
G(\gg) \bigr) \Bigr) G(\gg) - nG_s (\gg)
\end{eqnarray*}
Together with $f'\le 2$ and \eqref{ie:gve} we obtain
\[
\ca_{nG_s}(\gg) \,\le\, 4n G(\gg) \,\le\, 2n\gve^2 \,<\, 
\frac n \gs \,\le\, a(1) \,\le\, a(s) .
\]

\s Assume finally that $\gg$ lies in $D(2) \setminus D(\gve)$.
Then $\tau(|p(\gg)|) = 0$ and $f \circ G (\gg) = G(\gg)$.
Hence $G_-(\gg) = G(\gg)$ and
\begin{eqnarray*}
nG_s(\gg) &=& n \Bigl( \bigl(1-\gb(s) \bigr) G
+ \gb(s) \gs G_- \Bigr) (\gg) \\
&=& \bigl( 1+\gb(s)(\gs -1) \bigr) \2 nG (\gg).
\end{eqnarray*}
If $\ca_{nG_s}(\gg) = a(s)$, then $a=a(0) \in \cs(nG,q_0,q_1)$ in
view of Lemma~\ref{l:action}\:(ii) and the definition~\eqref{d:fas}
of $a(s)$, contradicting our hypothesis.
\proofend

\subsection{Almost all fibers are good} 
Fix now $q_0 \in M$. For $q \in M$ define
$$
D_q(4) \,=\, \left\{ p \in T_q^*M \mid |p| \le 4 \right\}
\,=\, D(4) \cap T_q^* M .
$$
Fix $n \in \NN$.
Given $H \in \ch_4(nG_+)$ let $U (q_0,H)$ be the set of those
$q_1 \in M$ for which $\gf_H (D_{q_0}(4) )$ and $D_{q_1}(4)$
intersect transversely.
Let $\mu_g$ be the Riemannian measure on $(M,g)$.
We say that a set $U \subset M$ has {\it full measure}\/
if $\mu_g(U) = \mu_g(M)$. Note that this property does not depend on
the Riemannian metric~$g$.

\begin{lemma}\label{l:full}
The set $U (q_0,H)$ is open and of full measure in~$M$.
\end{lemma}

\proof 
Since $D_{q_0}(4)$ is compact, $U(q_0,H)$ is open.
Applying Sard's Theorem to the projection $\gf_H (D_{q_0}(4)) \to
M$ one sees that $U(q_0,H)$ has full measure in~$M$.
\proofend

Note that for each $n \in \NN$ we have $\gf_{nH} = \gf^n_H$. 
Let $G_-$, $K$, $G_+$ be the Hamiltonians defined in~\eqref{def:G-},
\eqref{def:K}, \eqref{def:G+}.
Define
$$
V_n(q_0) \,:=\, U(q_0,nG_-) \cap U(q_0,nK) \cap U(q_0,nG_+).
$$
For $q_1 \in V_n(q_0)$ all the sets 
$\gf_{G_-}^n \bigl( D_{q_0}(4) \bigr)$, 
$\gf_K^n \bigl( D_{q_0}(4) \bigr)$, 
$\gf_{G_+}^n \bigl( D_{q_0}(4) \bigr)$ 
intersect $D_{q_1}(4)$ transversely. Moreover, 
\begin{corollary} \label{c:full}
The sets $V_n(q_0)$ have full measure in~$M$.
\end{corollary}

\section{Lagrangian Floer homology} \label{s:Floer}

\ni Floer homology for Lagrangian intersections was invented by
Floer in a series of seminal papers, \cite{F2,F3,F1,F4}. We shall
use a version of Lagrangian Floer homology described
in~\cite{KS,FS:GAFA,FS3}.

\subsection{Definition of $\HF_*^a(H, q_0,q_1;\FF_p)$}

\subsubsection{The chain groups}
Let $H \in \ch_4(nG_+)$, let $q_0 \in M$, and fix $q_1 \in U(q_0,H)$. 
For $a < n+1$ define
\[
\cp^a (H,q_0,q_1) \,:=\, \left\{ x \in \cp (H,q_0,q_1) \mid \ca_H
(x) < a \right\} .
\]
For $\gg \in \cp (H,q_0,q_1)$ outside $D(4)$ we have, by
Lemma~\ref{l:action}\:(i),
\[
\ca_H (\gg) \,=\, nG_+(\gg) \,=\, n \gs G(\gg) \,\ge\, 2n \,\ge\, n+1
\]
whence
\begin{equation}  \label{inc:Prho}
\cp^a (H,q_0,q_1) \,\subset\, D(4) .
\end{equation}
Since $q_1 \in U(q_0,H)$ we conclude that $\cp^a (H,q_0,q_1)$ is a
finite set.
The fibers of $T^*M$ form a Lagrangian foliation. 
For each path $x \in \cp(H,q_0,q_1)$ the Maslov index $\mu (x)$ is therefore a
well-defined integer. In case of a Riemannian Hamiltonian, $\mu (x)$
agrees with the Morse index of the corresponding geodesic path,
see~\cite[Section~1.2]{AS}.
Define the $k^{th}$ Floer chain group $\CF_k^a \left( H,
q_0,q_1;\FF_p \right)$ as the finite-dimensional $\FF_p$-vector
space freely generated by the elements of $\cp^a \left( H,q_0,q_1
\right)$ of Maslov index $k$, and define the full Floer chain group
as
\[
\CF_*^a \left( H, q_0,q_1;\FF_p \right) \,=\, \bigoplus_{k\in \ZZ}
\CF_k^a \left( H, q_0,q_1;\FF_p \right) .
\]

\subsubsection{Almost complex structures}
Recall that an almost complex structure $J$ on $T^*M$ is {\it
$\go$-compatible}\/ if
\[
\langle \cdot, \cdot \rangle \,\equiv\, g_J (\cdot, \cdot) :=
\go(\cdot, J\cdot)
\]
defines a Riemannian metric on $T^*M$. 
Denote by $\xi = \ker \left( \iota_Y\go |_{\pp D(4)}\right)$ the contact structure 
on $\pp D(4)$ defined by~$Y$ and $\go$, and for $t \ge 0$ let $\psi_t$ be
the semi-flow of the Liouville vector field $Y = \sum p_j \;\!\pp
p_j$ on $T^*M \setminus \Dcirc (4)$. 
An $\go$-compatible almost complex structure~$J$ on $T^*M$ is {\it convex}\/ 
on $T^*M \setminus \Dcirc(4)$ if
$$
J\xi = \xi, \qquad \go \bigl( Y(x),J(x)Y(x) \bigr) = 1, \qquad
D\psi_t (x) J(x) = J(\psi_t(x)) D\psi_t(x)
$$
for all $x \in \pp D(4)$ and $t \ge 0$. Following
\cite{CFH,V2,BPS} we consider the set $\cj$ of $t$-dependent smooth
families $\bJ = \{ J_t\}$, $t \in [0,1]$, of \text{$\go$-compatible}
almost complex structures on $T^*M$ such that $J_t$ is convex and
independent of $t$ on $T^*M \setminus \Dcirc(4)$. 
The set~$\cj$ is non-empty and connected.

\subsubsection{The boundary operators}
For $\bJ \in \cj$, for smooth maps $u$ from the strip $S = \RR
\times [0,1]$ to $T^*M$, and for $x^\pm \in \cp^a \left( H,q_0,q_1
\right)$ consider the partial differential equation with boundary
conditions
\begin{eqnarray}  \label{e:Floer}
 \left\{
  \begin{array}{lcr}
    \pp_s u-J_t(u) \left( \pp_t u - X_H(u) \right) = 0,  \\ [0.4em]
    u(s,j) \in T^*_{q_j}M,\,\, j=0,1 , \\ [0.4em]
    \displaystyle{\lim_{s \to \pm \infty}} u(s,t) \,=\, x^\pm (t)
             \;\text{ uniformly in } t .
  \end{array}
 \right.
\end{eqnarray}
The Cauchy-Riemann equation in~\eqref{e:Floer} is called {\it
Floer's equation}.

\begin{lemma}  \label{l:convex}
Solutions of equation~\eqref{e:Floer} are contained in $D(4)$.
\end{lemma}

\ni {\it Sketch of proof.} By~\eqref{inc:Prho} we have $x^\pm
\subset D(4)$, whence
\begin{equation}  \label{e:limit}
\lim_{s \to \pm \infty} u(s,t) \,=\, x^\pm (t) \subset D(4) .
\end{equation}
In view of the strong maximum principle, the lemma follows from the
convexity of $J$ outside $D(4)$ and from~\eqref{e:limit} together
with the fact that~$H = nG_+$ outside $D(4)$ implies $\go \left(
Y, J X_H \right) =0$, cf.~\cite{KS,FS:GAFA}. 
\proofend


We denote the set of solutions of~\eqref{e:Floer} by $\cm \left(
x^-,x^+,H;\bJ \right)$. Note that the group $\RR$ freely acts on
$\cm \left( x^-,x^+,H;\bJ \right)$ by time-shift.
Lemma~\ref{l:convex} is an important ingredient to establish the
compactness of the quotients $\cm \left( x^-,x^+,H;\bJ \right) /
\RR$. The other ingredient is that there is no bubbling-off of
$\bJ$-holomorphic spheres or discs. Indeed, $[\go]$ vanishes on
$\pi_2 \left( T^*M \right)$ because $\go = d \gl$ is exact, and
$[\go]$ vanishes on $\pi_2 \bigl( T^*M , T_{q_j}^*M \bigr)$ because
$\gl$ vanishes on $T_{q_j}^*M$, $j=0,1$. See for instance~\cite{F3}
or \cite{Sa-90} for more details.

There exists a residual subset $\cj_{\reg}(H)$ of $\cj$ such that for
each $\bJ \in \cj_{\reg}(H)$ the linearized operator for Floer's
equation is surjective for each solution of~\eqref{e:Floer}. For
such a {\it regular}\/ $\bJ$ the moduli space $\cm \left(
x^-,x^+,H;\bJ \right)$ is a smooth manifold of dimension $\mu(x^-) -
\mu(x^+)$ for all $x^\pm \in \cp^a \left( H,q_0,q_1 \right)$,
see~\cite{FHS}. Fix $\bJ \in \cj_{\reg}(H)$. It is shown
in~\cite[Section~1.4]{AS} that the manifolds $\cm \left(
x^-,x^+,H;\bJ \right)$ can be oriented in a way which is coherent
with gluing. For $x^\pm \in \cp^a \left( H, q_0,q_1 \right)$ with
$\mu(x^-) = \mu(x^+)+1$ let
\[
n(x^-,x^+,H;\bJ) \,=\, \# \cm (x^-, x^+,H;\bJ) / \RR \,\in\, \ZZ
\]
be the oriented count of the finite set $\cm (x^-, x^+,H;\bJ) /
\RR$. If $u \in \cm \left( x^-,x^+,H;\bJ \right)$, then $\ca_H(x^-)
\ge \ca_H(x^+)$, see the more general Lemma~\ref{l:homotopy} below.
For $k \in \ZZ$ one can therefore define the Floer boundary operator
\[
\pp_k \left( \bJ \right) \colon \CF_k^a \left( H, q_0,q_1;\FF_p
\right) \,\to\, \CF_{k-1}^a \left( H, q_0,q_1;\FF_p \right)
\]
as the linear extension of
\[
\pp_k(\bJ) \, x^- \,=\, \sum n(x^-,x^+,H;\bJ) \, x^+
\]
where $x^- \in \cp^a(H, q_0,q_1)$ has index $\mu(x^-)=k$
and the sum runs over all $x^+ \in \cp^a( H, q_0,q_1)$ of index
$\mu(x^+)=k-1$.
Then $\pp_{k-1} \left( \bJ \right) \circ \pp_k \left( \bJ \right)
=0$ for each $k$. The proof makes use of the compactness of the $0$-
and $1$-dimensional components of $\cm \left( x^-,x^+,H; \bJ \right)
/ \RR$, see~\cite{F2,Schw,AS}.
As our notation suggests, the Floer homology groups
\[
\HF_k^a \left( H, q_0,q_1;\FF_p \right) \,:=\, \ker \pp_k (\bJ) / \im
\pp_{k+1} (\bJ)
\]
do not depend on the choices involved in their construction: They
neither depend on coherent orientations up to canonical
isomorphisms, \cite[Section~1.7]{AS}, nor do they depend on $\bJ \in
\cj_{\reg}(H)$ up to natural isomorphisms, as a continuation
argument shows, \cite{F2,Schw}. The groups $\HF^a_k( H, q_0,q_1;\FF_p)$
do depend, however, on $a < n+1$ and $H \in \ch_4(nG_+)$.
In the sequel, the points $q_0,q_1 \in M$ and the field~$\FF_p$ are
fixed throughout. We shall therefore often write $\CF_*^a (H)$ and
$\HF_*^a (H)$ instead of $\CF_*^a ( H, q_0,q_1;\FF_p)$ and $\HF_*^a
(H, q_0,q_1;\FF_p)$.

\subsection{Continuation homomorphisms}
The goal of this section is to relate the groups
$$
\HF_*^a (nG_-), \quad \HF_*^a (nK), \quad \HF_*^a (nG_+) .
$$
Let $\gb \colon \RR \to [0,1]$ be the function
from~\eqref{e:beta}. Given two functions $H^-, H^+ \in
\ch_4(nG_+)$ with $H^-(x) \le H^+(x)$ for all $x \in T^*M$ we
form the monotone homotopy
\begin{equation}  \label{e:Hs}
H_s (x) \,=\, H^-(x) +\gb(s) \bigl( H^+(x)-H^-(x) \bigr) .
\end{equation}
Then $H_s \in \ch_4(nG_+)$ for each $s$, and $H_s = H^-$ for $s
\le 0$ and $H_s = H^+$ for $s \ge 1$.
Consider the equation with boundary conditions
\begin{eqnarray}  \label{e:Floers}
 \left\{
  \begin{array}{lcr}
    \pp_s u-J_{s,t}(u) \left( \pp_t u - X_{H_s}(u) \right) = 0,  \\ [0.4em]
    u(s,j) \in T^*_{q_j}M,\,\, j=0,1 , \\ [0.4em]
    \displaystyle{\lim_{s \to \pm \infty}} u(s,t) \,=\, x^\pm (t)
             \;\text{ uniformly in } t .
  \end{array}
 \right.
\end{eqnarray}
where $s \mapsto \{ J_{s,t} \}$, $s \in \RR, t \in [0,1]$, is a {\it
regular homotopy}\/ of families $\{J_t \}$ of almost complex
structures on $T^*M$. This means that
\begin{itemize}
\item[$\bullet$]
$J_{s,t}$ is $\go$-compatible and convex and independent of $s$ and
$t$ outside $D(4)$;
\item[$\bullet$]
$J_{s,t} = J_t^- \in \cj_{\reg} (H^-)$ for $s \le 0$;
\item[$\bullet$]
$J_{s,t} = J_t^+ \in \cj_{\reg} (H^+)$ for $s \ge 1$;
\item[$\bullet$]
the solutions of~\eqref{e:Floers} are transverse (that is, the
associated Fredholm operators are surjective) and therefore form
finite dimensional moduli spaces.
\end{itemize}
The following lemma is well-known. We reprove it in view of the many
different sign conventions 
(for the Hamilton equation, the action functional, the Floer equation, etc.) 
used by different authors.

\begin{lemma}  \label{l:homotopy}
Assume that $u \colon S = \RR \times [0,1] \to T^*M$ is a solution
of~\eqref{e:Floers}. Then
\[
\ca_{H^+}(x^+) \,\le\, \ca_{H^-}(x^-) -
\int_0^1\int_{-\infty}^\infty \gb'(s) \left( H_1-H_0 \right)(u(s,t))
\, ds\,dt .
\]
\end{lemma}

\proof 
Note that $\gl = p \2 dq$ vanishes along the Lagrangian boundary conditions
$u(s,j) \in T_{q_j}^*M$ in~\eqref{e:Floers}.
Since $\go = d\gl$ is exact, we obtain, using Stokes Theorem and taking into
account orientations,
$$
\int_S u^*\go \,=\, \int_{u (\pp S)} \gl
\,=\, \int_{x^+} \gl - \int_{x^-} \gl .
$$
Moreover, by~\eqref{e:Hs},
\[
\frac{d}{ds}H_s (u(s,t)) \,=\, dH_s (u(s,t)) (\pp_su) + \gb'(s)
(H_1-H_0) (u(s,t)).
\]
This and the asymptotic boundary condition in~\eqref{e:Floers} yield
\begin{eqnarray*}
\int_0^1 \int_{-\infty}^\infty  dH_s (u(s,t)) (\pp_su)\,ds\,dt &=&
\int_0^1 H^+(x^+) \,dt - \int_0^1 H^-(x^-) \,dt \\
&&- \int_0^1 \int_{-\infty}^\infty \gb'(s) (H_1-H_0)(u(s,t))
\,ds\,dt.
\end{eqnarray*}
Together with the compatibility $g_{s,t}(v,w) = \go(v,J_{s,t}w)$,
Floer's equation in~\eqref{e:Floers}, Hamilton's
equation~\eqref{e:XH} and the definition~\eqref{def:af} of the
action functional we obtain
\begin{eqnarray*}
0 &\le&
\int_0^1 \int_{-\infty}^\infty g_{s,t} (\pp_s u, \pp_s u) \,ds\,dt \\
&=& \int_0^1 \int_{-\infty}^\infty g_{s,t} \bigl(\pp_su, J_{s,t}(u)
\left( \pp_tu - X_{H_s}(u)\right) \bigr) \,ds\,dt\\
&=& - \int_S u^* \go - \int_0^1 \int_{-\infty}^\infty \go \left(
X_{H_s}(u), \pp_su \right) ds\, dt \\
&=& \int_{x^-} \gl - \int_{x^+} \gl + \int_0^1 \int_{-\infty}^\infty
dH_s(u) (\pp_su) \,ds\, dt \\
&=& \int_{x^-} \gl - \int_{x^+} \gl + \int_0^1 H^+(x^+) \,dt -
\int_0^1 H^-(x^-) \,dt  \\
&& \qquad \qquad \qquad \, - \int_0^1 \int_{-\infty}^\infty \gb'(s)
(H_1-H_0)(u(s,t)) \,ds\,dt \\
&=& \ca_{H^-}(x^-) - \ca_{H^+}(x^+) -\int_0^1 \int_{-\infty}^\infty
\gb'(s) (H_1-H_0)(u(s,t)) \,ds\,dt,
\end{eqnarray*}
as claimed. 
\proofend

In view of Lemma~\ref{l:homotopy} the action decreases along
solutions $u$ of~\eqref{e:Floers}. By counting these solutions one
can therefore define the Floer chain map
$$
\phi_{H^+H^-} \colon \CF^a_* (H^-) \to \CF^a_* (H^+),
$$
see~\cite{F5,SZ,FH,CFH}. The induced {\it continuation homomorphism}
$$
\Phi_{H^+H^-} \colon \HF^a_* (H^-) \to \HF^a_* (H^+)
$$
on Floer homology does not depend on the choice of the regular
homotopy $\{J_{s,t}\}$ used in its definition. An important property
of these homomorphisms is naturality with respect to concatenation,
\begin{equation}  \label{e:concatenation}
\Phi_{H_3H_2} \circ \Phi_{H_2 H_1} \,=\, \Phi_{H_3H_1} \quad
\text{for }\, H_1 \le H_2 \le H_3 .
\end{equation}
Another important fact is the following invariance property, which
is proved in~\cite{FH,CFH} and \cite[Section~4.5]{BPS}.
\begin{lemma}  \label{l:invariance}
If $a \notin \cs (H_s,q_0,q_1)$ for all $s \in [0,1]$, then
$\Phi_{H^+H^-} \colon \HF_*^a(H^-) \to \HF_*^a(H^+)$ is an
isomorphism.
\end{lemma}

After these recollections we return to our Hamiltonians 
$nG_- \le nK \le nG_+$. 
Fix $q_1 \in V_n(q_0)$. 
By~\eqref{inc:Prho}, the set
$$
\cs (n) \,:=\, \bigl( \cs (nG_-,q_0,q_1) \cup \cs (nK,q_0,q_1)
\bigr) \cap (n,n+1)
$$
is finite. In particular,
\[
\gd_n \,:=\, \min \left\{ s \in \cs (nK,q_0,q_1) \mid s > n \right\}
\,>\, 0 .
\]
Choose
\begin{equation}  \label{def:am}
a_n \,\in\, \bigl( (n,n+1) \cap (n,n+\gd_n) \bigr) \setminus \cs (n) .
\end{equation}
Then $a_n$ is not in the action spectrum of $nG_-$, $nK$, and by
Proposition~\ref{p:action},
\begin{equation}  \label{e:iff}
x \in \cp(nK,q_0,q_1) \cap D(\Sigma) \quad \Longleftrightarrow \quad
\ca_{nK} (x) < a_n .
\end{equation}
%


Our next goal is to show that $\HF^{a_n}_* (nG_-)$ and
$\HF^{a_n/\gs}_*(nG_+)$ are naturally isomorphic. This is a special
case of a generalization of Lemma~\ref{l:invariance} stated as
Proposition~1.1 in~\cite{V2}. We give an ad hoc construction in the
case at hand.

\begin{figure}[h]
 \begin{center}
  \psfrag{s}{$s$}
  \psfrag{t}{$t$}
  \psfrag{an}{$a_n$}
  \psfrag{bn}{$b_n$}
  \psfrag{fa}{$\frak{a}_1$}
  \psfrag{ans}{$\frac{a_n}{\gs}$}
  \psfrag{s0}{$s_0=0$}
  \psfrag{s1}{$s_1$}
  \psfrag{sk}{$s_k$}
  \psfrag{sk1}{$s_{k+1}=1$}
  \leavevmode\epsfbox{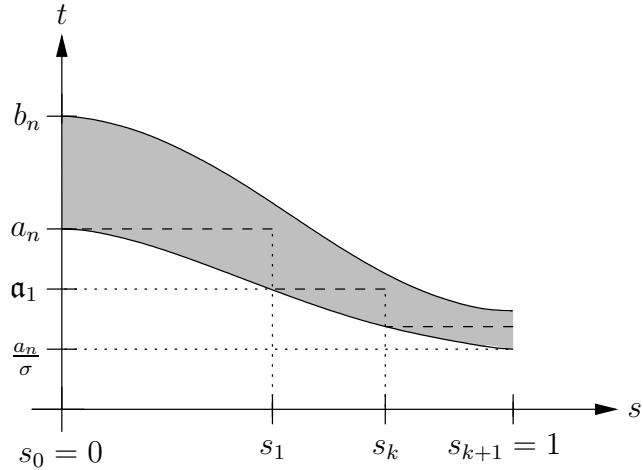}
 \end{center}
 \caption{The curves $a_n(s)$ and $b_n(s)$.}
 \label{figure.ambm}
\end{figure}

Choose $b_n \in (a_n,n+1)$ such that
\begin{equation}  \label{e:amhat}
[a_n,b_n] \cap \cs (nG_-) \,=\, \emptyset ,
\end{equation}
and for $s \in [0,1]$ define
$$
a_n(s) = \frac{a_n}{1+\gb(s)(\gs -1)},  \qquad \ b_n(s) =
\frac{b_n}{1+\gb(s)(\gs -1)} .
$$
For $(s,t)$ in the gray band bounded by the graphs of $a_n(s)$ and
$b_n(s)$ we have $t \notin \cs(nG_s)$ in view of the Non-crossing
Lemma~\ref{l:NC}.
Choose a partition $0 = s_0 < s_1 < \dots < s_k < s_{k+1} = 1$ so fine that
\begin{equation}  \label{ie:amhat}
b_n(s_{j+1}) \,>\, a_n(s_j), \qquad j=0,\dots,k.
\end{equation}
Abbreviate $\frak{a}_j = a_n(s_j)$ and $\frak G_j = n G_{s_j}$. Then
\begin{equation}  \label{e:amnot}
\frak a_j \notin \cs (nG_s) \quad \text{for }\, s \in [s_j,s_{j+1}]
,
\end{equation}
cf.~Figure~\ref{figure.ambm}.

%
%

\ni
Together with Lemma~\ref{l:invariance}
we find that
\[
\Phi_{\frak G_{j+1}\frak G_j} \colon \HF_*^{\frak a_j} (\frak G_j)
\to \HF_*^{\frak a_j}(\frak G_{j+1})
\]
is an isomorphism. Since $[\frak a_{j+1},\frak a_j] \cap \cs (\frak
G_{j+1}) = \emptyset$, we have $\HF_*^{\frak a_j} (\frak G_{j+1}) =
\HF_*^{\frak a_{j+1}}(\frak G_{j+1})$, and so
\[
\widehat \Phi_{\frak G_{j+1}\frak G_j} \equiv \Phi_{\frak G_{j+1}
\frak G_j} \colon \HF_*^{\frak a_j} (\frak G_j) \to \HF_*^{\frak
a_{j+1}}(\frak G_{j+1})
\]
is an isomorphism.
Recalling $\frak a_0 = a_n$ and $\frak a_{k+1} = a_n/\gs$ we obtain
that the composition
\[
\widehat \Phi_{nG_+nG_-} \,:=\, \widehat \Phi_{\frak G_{k+1} \frak
G_k} \circ \dots \circ \widehat \Phi_{\frak G_2 \frak G_1} \circ
\widehat \Phi_{\frak G_1 \frak G_0} \colon \HF_*^{a_n}(nG_-) \to
\HF_*^{a_n/\gs}(nG_+)
\]
is an isomorphism.
Let
$$
\HF_*(\iota) \colon \HF_*^{a_n/\gs}(nG_+) \to \HF_*^{a_n}(nG_+)
$$
be the homomorphism induced by the inclusion $\CF_*^{a_n/\gs}(nG_+)
\to \CF_*^{a_n}(nG_+)$.

\begin{proposition}  \label{p:diagram}
For each $k$ and $n$ there is a commutative diagram of homomorphisms
$$  
\xymatrix{ 
&\HF_k^{a_n/\gs}(nG_+) \ar[dr]^{\HF_k(\iota)}& \\
\HF_k^{a_n}(nG_-) \ar[ur]^{\widehat{\Phi}_{nG_+nG_-}}
\ar[rr]^{\Phi_{nG_+nG_-}}
\ar[dr]_{\Phi_{nKnG_-}} && \HF_k^{a_n}(nG_+) \\
& \HF_k^{a_n}(nK) \ar[ur]_{\Phi_{nG_+nK}}&}
$$
and $\widehat \Phi_{nG_+nG_-}$ is an isomorphism.
\end{proposition}

\proof 
By construction, the isomorphism $\widehat \Phi_{nG_+nG_-}$
is induced by the composition of Floer chain maps
$$
\widehat \phi_{\frak G_{j+1} \frak G_j} \colon \CF_*^{\frak
a_j}(\frak G_j) \to \CF_*^{\frak a_{j+1}}(\frak G_{j+1}) \subset
\CF_*^{\frak a_j}(\frak G_{j+1}).
$$
Therefore, $\HF_*(\iota) \circ \widehat \Phi_{nG_+ nG_-}$ is induced by
the composition of Floer chain maps
$$
\phi_{\frak G_{j+1} \frak G_j} \colon \CF_*^{\frak a_0}(\frak G_j)
\to \CF_*^{\frak a_0}(\frak G_{j+1}) .
$$
By~\eqref{e:concatenation}, this composition induces the same map in
Floer homology as
$$
\phi_{nG_+nG_-} \colon \CF_*^{\frak a_0}(nG_-) \to \CF_*^{\frak
a_0}(nG_+).
$$
The upper triangle therefore commutes. The lower triangle commutes
in view of $nG_- \le nK \le nG_+$ and according
to~\eqref{e:concatenation}. 
\proofend

\begin{corollary} \label{c:floer}
$\dim \HF_k^{a_n} (nK,q_0,q_1;\FF_p) \,\ge\, \rank \HF_k(\iota) \colon 
\HF_k^{a_n/\gs}(nG_+) \to \HF_k^{a_n}(nG_+)$.
\end{corollary}

\subsection{From Floer homology to the homology of the based loop space}
\label{s:Floer.based}

In this section we use Corollary~\ref{c:floer}
and our assumption that~$M$ is energy hyperbolic to prove

\begin{theorem}  \label{t:floer}
Let $(M,g)$ and $K$ be as above.
Then there exist $h>0$, $p \in \PP$ and $N \in \NN$ depending only on $(M,g)$ 
such that for all $n \ge N$ and for $a_n$ as above the following holds.
For all $q_0 \in M$ and all $q_1 \in V_n(q_0)$,
$$
\dim \HF_*^{a_n} \left( nK, q_0, q_1 ;\FF_p \right)
\,\ge\,  e^{h \1 n} .
$$
\end{theorem}

\proof
For $q_0, q_1 \in M$ let $\Omega^1 \left( M,q_0,q_1\right)$
be the space of all paths $q \colon [0,1] \to M$ of Sobolev
class $W^{1,2}$ such that $q(0)=q_0$ and $q(1)=q_1$.
Again, this space has a canonical Hilbert manifold structure.
The energy functional $\ce \colon \Omega^1 \left( M,q_0,q_1\right) \to \RR$
is defined as
$$
\ce (q) \,=\, \frac{1}{2} \int_0^1 \left| \dot{q}(t) \right|^2 \,dt.
$$
For $b>0$ we consider the sublevel sets $\ce^b(q_0,q_1) := \left\{ q
\in \Omega^1 (M,q_0,q_1 ) \mid \ce(q) \le b \right\}$.
\begin{proposition} \label{p:AS}
For each $k$ and $n$ there is a commutative diagram of homomorphisms
$$
\xymatrix{ 
\HF_k^{a_n/\gs} (nG_+,q_0,q_1;\FF_p)
\ar[d]_{\HF_k(\iota)} \ar[r]^-{\cong} &
\H_k(\ce^{na_n}(q_0,q_1);\FF_p) \ar[d]^{\H_k(\iota)} \\
\HF_k^{a_n} (nG_+,q_0,q_1;\FF_p) \ar[r]^-{\cong} & 
\H_k(\ce^{\gs na_n}(q_0,q_1);\FF_p) 
}
$$
where the horizontal maps are isomorphisms and the right map 
$\H_k(\iota)$ is induced by the inclusion $\ce^{na_n}(q_0,q_1)
\ha \ce^{\gs n a_n}(q_0,q_1)$.
\end{proposition}

\proof 
Let $L \colon TM \to \RR$ be the Legendre transform of
$nG_+$, let
\[
\ce_L (q) \,=\, \int_0^1 L \bigl( q(t), \dot{q} (t) \bigr) dt
\]
be the corresponding action functional on $\Omega^1 (M,q_0,q_1)$,
and let
\[
\ce^b_L(q_0,q_1) \,=\, \left\{ q \in \Omega^1 (M,q_0,q_1) \mid
\ce_L(q) \le b \right\} .
\]
Applying the Abbondandolo--Schwarz Theorem~3.1 of~\cite{AS} to
$nG_+$ and $L$, we obtain for each $b>0$ an isomorphism
$$
\xymatrix{ \HF_k^b \left(nG_+,q_0,q_1;\FF_p \right) & \HM_k^b
(L,q_0,q_1;\FF_p) \ar[l]_-{\;\;\;\Theta_k^b}}.
$$
Here, $\HM_k^b (L,q_0,q_1;\FF_p)$ is the Morse homology ``below
level $b$" of $\ce_L$ constructed in~\cite{AM}, see also
\cite[Section~2]{AS}. The Abbondandolo--Schwarz chain isomorphisms
\[
\xymatrix{ \CF_k^b \left(nG_+,q_0,q_1;\FF_p \right) & \CM_k^b
(L,q_0,q_1;\FF_p) \ar[l]_-{\;\;\;\theta_k^b}}
\]
between the Morse and the Floer chain complexes commute with
inclusions
$$
\CF_k^b(nG_+,q_0,q_1;\FF_p) \ha \CF_k^{b'} (nG_+,q_0,q_1;\FF_p)
\quad \text{and} \quad \CM_k^b(L,q_0,q_1;\FF_p) \ha
\CM_k^{b'}(L,q_0,q_1;\FF_p)
$$
for $b < b'$, see~\cite[p.~298]{AS}. Therefore, the induced diagram
of homology groups commutes,
$$
\xymatrix{ \HF_k^{b}(nG_+,q_0,q_1;\FF_p) \ar[d]_{\HF_k(\iota)} &
\HM_k^b(L,q_0,q_1;\FF_p) \ar[d]^{\HM_k(\iota)} \ar[l]_-{\;\;\Theta_k^b}\\
\HF_k^{b'}(nG_+,q_0,q_1;\FF_p) & \HM_k^{b'}(L,q_0,q_1;\FF_p)
\ar[l]_-{\;\;\Theta_k^{b'}} }
$$
Moreover, Abbondandolo and Majer constructed chain isomorphisms
$$
\xymatrix{ \CM_k^b \left(L,q_0,q_1;\FF_p \right) \ar[r]^-{\tau_k^b}
& \C_k (\ce^b_L(q_0,q_1);\FF_p) }
$$
between the singular chain complexes and the Morse chain complexes,
which commute with inclusions
$$
\CM_k^b(L,q_0,q_1;\FF_p) \ha \CM_k^{b'}(L,q_0,q_1;\FF_p) \quad
\text{and} \quad \C_k(\ce_L^b(q_0,q_1);\FF_p) \ha \C_k
(\ce_L^{b'}(q_0,q_1);\FF_p)
$$
for $b < b'$, see~\cite{AM} and \cite[Section~2.3]{AS}. Therefore,
the induced diagram of homology groups commutes,
$$
\xymatrix{ \HM_k^{b}(L,q_0,q_1;\FF_p) \ar[r]^{T_k^b}
\ar[d]_{\HM_k(\iota)} &
\H_k(\ce^b_L(q_0,q_1);\FF_p) \ar[d]^{\H_k(\iota)} \\
\HM_k^{b'}(L,q_0,q_1;\FF_p) \ar[r]^{T_k^{b'}} &
\H_k(\ce^{b'}_L(q_0,q_1);\FF_p)}
$$
Notice now that $L(q,v) = \frac{1}{n \gs} \frac 12
\left|v\right|^2$, whence $\ce^b_L (q_0,q_1) = 
\ce^{n \gs b}(q_0,q_1)$ for all $b>0$. The proposition follows. 
\proofend

Consider now the commutative diagram
$$
\xymatrix{ 
\H_k \bigl(\ce^{na_n}(q_0,q_1);\FF_p \bigr) 
\ar[d] \ar[dr]^{\iota_k} & \\
\H_k \bigl(\ce^{\gs na_n}(q_0,q_1);\FF_p \bigr) 
\ar[r] &
\H_k \bigl( \Omega^1(M,q_0,q_1);\FF_p \bigr)
}
$$                     
induced by the inclusions 
$\ce^{na_n}(q_0,q_1) \subset \ce^{\gs na_n}(q_0,q_1) \subset \Omega^1(M,q_0,q_1)$.
%

\begin{lemma}  \label{l:homeq}
For each $c \in (0,1)$ there exists $A>0$ depending only on $c$ and $(M,g)$ such that
$$
\dim \iota_k \H_k \left( \ce^a (q_0,q_1); \FF_p \right) \,\ge\,
\dim \iota_k \H_k \left( \ce^{ca} (q_0); \FF_p \right)
$$
for all $k$ and all $a \ge A$.
\end{lemma}

\proof
Let $\rho$ be the diameter of $(M,g)$.
Choose a path $\frak p  \colon [0,1] \to M$ from $q_0$ to~$q_1$ of length $\le \rho$.
Parametrize $\frak p$ proportional to arc-length.
Then $\ce (\frak p) \le \frac 12 \rho^2$.
Let $\tau := \frac 12 (c+1) \in (0,1)$.
For $\gg \in \Omega^1(M,q_0)$ define $\gg \ast_\tau \frak p \in \Omega^1(M,q_0,q_1)$ by
\begin{eqnarray*}  
\left( \gg \ast_\tau \frak p \right) (t) \,=\,
 \left\{       
  \begin{array}{ll}
    \gg \left( \tfrac t \tau \right), & 0 \le t \le \tau,  \\ [0.2em]
    \frak p \left( \tfrac{t-\tau}{1- \tau} \right), & \tau \le t \le 1. 
 \end{array}
 \right.
\end{eqnarray*}
The map $\frak P \colon \Omega^1(M,q_0) \to \Omega^1(M,q_0,q_1)$, $\gg \mapsto \gg \ast_\tau \frak p$, 
is a homotopy equivalence with homotopy inverse 
$\Omega^1(M,q_0,q_1) \to \Omega^1(M,q_0)$, $\delta \mapsto \delta \ast_\tau \frak p^{-1}$.
Notice that 
\[
\ce \left( \gg \ast_\tau \frak p \right) \,=\, 
\tfrac 1 \tau \ce (\gg) + \tfrac{1}{1-\tau} \ce (\frak p) 
\quad\, \text{ for all }\, \gg \in \Omega^1(M,q_0) .
\]
Therefore, 
$$
\frak P \left( \ce^{ca} (q_0) \right) \,\subset\, 
\ce^{ \frac{ca}{\tau} + \frac{\ce(\frak p)}{1-\tau} } (q_0,q_1) 
\quad\, \text{ for all }\, a .
$$
Since $\ce (\frak p) \le \frac 12 \rho^2$, there exists $A>0$ depending only on
$c$ and $\rho$ such that
$$
\frac{ca}{\tau} + \frac{\ce(\frak p)}{1-\tau} \,\le\, a  
\quad\, \text{ for all }\, a \ge A .
$$
Hence 
$$
\frak P \left( \ce^{ca} (q_0) \right) \,\subset\, 
\ce^{a} (q_0,q_1) 
\quad\, \text{ for all }\, a \ge A .
$$
Since $\frak P$ is a homotopy equivalence,
it follows that
$$
\dim \iota_k \H_k \left( \ce^{ca}(q_0); \FF_p \right) 
\,\le\,
\dim \iota_k \H_k \left( \ce^a (q_0,q_1); \FF_p \right) 
\quad\, \text{ for all }\, a \ge A ,
$$
as claimed.
\proofend

Since $(M,g)$ is energy hyperbolic, $h := C(M,g) >0$.
By definition of $C(M,g)$ and by Lemma~\ref{l:homeq}
there exist $p \in \PP$ and $N_0 \in \NN$
such that for all $m \ge N_0$, 
$$
\sum_{k \ge 0} \dim \iota_k \H_k \bigl( \ce^{\frac 12 m^2} (q_0,q_1) ;\FF_p \bigr) 
\,\ge\,
e^{\frac 1{\sqrt 2} \1 h \2 m} .
$$
Recall that $a_n \ge n$.
Therefore there exists $N \in \NN$ such that for all $n \ge N$,
$$
\sum_{k \ge 0} \rank \iota_k \,\ge\, e^{h \1 n} . 
$$
Together with Proposition~\ref{p:AS} we find that
$$
\sum_{k \ge 0} \rank \HF_k(\iota) \,\ge\, \sum_{k \ge 0} \rank \iota_k \,\ge\, e^{h \1 n} . 
$$
This and Corollary~\ref{c:floer} yield Theorem~\ref{t:floer}.
\proofend

\begin{remark*}
{\rm
The essential point in the proof of Theorem~\ref{t:floer} 
is that the sum $\sum_{k \ge 0} \rank \iota_k$ grows exponentially with~$n$.
This is so by our assumption that~$M$ is energy hyperbolic.
In the special case that the fundamental group of~$M$ has exponential growth,
$\rank \iota_0$ grows exponentially with~$n$.
In another special case where~$M$ is simply connected and hyperbolic,
a theorem of Gromov~\cite{Gromov78} 
guarantees the existence of a constant~$c(M,g)$
such that $\iota_k$ is surjective if $n \ge c(M,g) k$,
whence $\sum_{k \ge 0} \rank \iota_k$ grows exponentially with~$n$.
}
\end{remark*}

\section{Proof of Theorem 1} \label{s:th1}

\subsection{From the growth of Floer homology to volume growth of $D_{q_0}(\Sigma)$}
\label{s:Floer.growth}

Theorem~\ref{t:floer} implies that the Riemannian volume of the
sequence of submanifolds $\gf_K^n \bigl( D_{q_0}(\Sigma)\bigr)$ grows exponentially.
Indeed, 
fix $q_0 \in M$, and let $p \in \PP$ and $N \in \NN$ be as in Theorem~\ref{t:floer}. 
Let $n \ge N$ and 
pick $q_1 \in V_n (q_0)$. 
By~\eqref{e:iff},
the generators of $\CF_k^{a_n} \left( nK,q_0,q_1 ;\FF_p \right)$ 
correspond to
$\gf_K^n \left( D_{q_0}(\Sigma) \right) \cap D_{q_1}(\Sigma)$.
Therefore, 
\begin{eqnarray}  \label{e:nC}
\# \bigl( \gf_K^n (D_{q_0}(\Sigma)) \cap D_{q_1}(\Sigma) \bigr) 
&=&
\dim \CF_*^{a_n} \left( nK,q_0,q_1;\FF_p \right)  \\
&\ge&
\dim \HF_*^{a_n} \left( nK,q_0,q_1; \FF_p \right) \notag \\
&\ge&
e^{h \1 n} . \notag
\end{eqnarray}
Let $\mu_g$ be the Riemannian measure on $(M,g)$.
Let~$g^*$ be the Riemannian metric on $T^*M$ induced by~$g$,
and for a submanifold $\gs \subset T^*M$
let $\mu_{g^*}(\gs)$ be the measure of~$\gs$ with respect to
the Riemannian measure on~$\gs$ given by the Riemannian metric
on~$\gs$ induced by $g^*$. 
Note that $\pi \colon T^*M \to M$ is a Riemannian submersion
with respect to the Riemannian metrics $g^*$ and $g$,
and recall from Corollary~\ref{c:full} that $V_n(q_0)$ is an open subset of~$M$
of full measure.
Therefore,
\begin{equation}  \label{e1}
\mu_{g^*} \bigl( \gf_K^n (D_{q_0}(\Sigma)) \bigr) \,\ge\, 
e^{h \1 n} \mu_g \bigl( V_n(q_0) \bigr) \,=\, \mu_g (M) \,e^{h \1 n} .
\end{equation}

\subsection{From volume growth of the sublevel to topological entropy at the level}  
\label{s:vol.top}

If $\gf$ is a $C^\infty$-smooth map of a compact manifold~$P$,
a geometric way of defining the topological entropy of~$\gf$ 
was found by Yomdin and Newhouse in their seminal works~\cite{Y} and \cite{Newh}:
Fix a Riemannian metric~$g$ on~$P$.
For $j \in \left\{ 1, \dots, \dim P \right\}$
denote by $S_j$ the set of smooth compact
(not necessarily closed) \text{$j$-dimensional} submanifolds of~$P$.
The {\it $j$'th volume growth}\, of $\gf$ is defined as
\[
v_j (\gf) \,=\,
 \sup_{\gs \in S_j} \liminf_{n \ra \infty} \frac 1n
    \log \mu_g \left( \gf^n ( \gs ) \right) ,
\]
and the {\it volume growth}\, of $\gf$ is defined as
\begin{equation} \label{def:v}
v (\gf) \,=\, \max_{1 \le j \le \dim P} v_j(\gf) .
\end{equation}
They do not depend on the choice of the Riemannian metric~$g$ used 
in their definition.
Yomdin proved in \cite{Y} that $h_{\top} (\gf) \ge v (\gf)$,
and Newhouse proved in~\cite{Newh} that 
$h_{\text{top}} (\gf) \le v (\gf)$.
Thus 
\begin{equation}  \label{e:yomdin}
h_{\top} (\gf) \,=\, v (\gf) .
\end{equation}
With $P = D(\Sigma)$ and $\gf = \gf_K |_{D(\Sigma)}$
we have, by~\eqref{e1}, that
$v_d \bigl( \gf_K |_{D(\Sigma)} \bigr) >0$.
Hence, by~\eqref{def:v} and \eqref{e:yomdin},
\begin{equation} \label{e:GD}
h_{\top} (\gf_K |_{D(\Sigma)}) >0 .
\end{equation}
Recall that $\Sigma = \pp D(\Sigma)$.
The topological entropy of a flow is defined as the topological entropy of its time~1~map.
Theorem~1 follows from~\eqref{e:GD} and the following

\begin{proposition} \label{p:atlevel}
$h_{\top} (\gf_K |_{\Sigma}) = h_{\top} (\gf_K |_{D(\Sigma)})$.
\end{proposition}

\proof
For $s \in [0,1]$ define the diffeomorphism $\psi_s$ of $T^*M$ by
$\psi_s (q,p) = (q,sp)$.
Abbreviate $s\1\Sigma = \psi_s(\Sigma)$.
Recall from~\eqref{e:fSigma} that $K = f \circ F$ on $D(\Sigma)$.
Also recall from~\eqref{d:f} that $f'(r)=1$ at $r=1$.
For $(q,p) \in \Sigma$ we thus have
\begin{equation} \label{e:timechange}
X_{f\circ F} (q,sp) \,=\, \bigl( f'(F(q,sp)) s \bigr) \, d\psi_s (q,p) X_{f\circ F}(q,p) 
\,=\, \gs(s) \,d\psi_s (q,p) X_{f\circ F}(q,p) .
\end{equation}
For the latter equality we recalled that $F(q,sp) = s^2F(q,p) = s^2$
and abbreviated 
$$
\gs(s) \,=\, f'(s^2)s .
$$
If $s \le \gve$, then $\gs (s) =0$, whence $X_{f \circ F} |_{s \1 \Sigma} =0$ and
$h_{\top} \left( \gf_{f \circ F} |_{s\1 \Sigma}\right) =0$.
If $s > \gve$, then $\gs (s) >0$, and the identity~\eqref{e:timechange}
shows that $\psi_s$ conjugates the flows 
$\gf^t_{f\circ F}|_{\Sigma}$ and $\gf^{\frac{1}{\gs(s)}t}_{f\circ F}|_{s \1 \Sigma}$.
Topological entropy is an invariant of conjugacy, and scales like 
$h_{\top}(\phi^{ct}) = |c| \2 h_{\top}(\phi^t)$ for $c \in \RR$.
Therefore, 
$$
h_{\top} \left( \gf^t_{f\circ F}|_{s \1 \Sigma} \right) \,=\,
\gs(s)\, h_{\top} \left( \gf^{\frac{1}{\gs(s)}t}_{f\circ F}|_{s \1 \Sigma} \right) \,=\,
\gs(s)\, h_{\top} \left( \gf^t_{f\circ F}|_{\Sigma} \right) .
$$
If $s^2 \ge \gve^2$, then $f'(s^2) =1$, whence $\gs(s)=s\le1$.
If $s^2 \in (\gve^2,\eps)$, then $f'(s^2) \le 2$, whence $\gs(s) < 2s < 2 \sqrt \gve < 1$.
If follows that
$$
\sup_{s \in [0,1]} h_{\top} \left( \gf_{f \circ F} |_{s \1 \Sigma} \right) \,=\, 
h_{\top} \left( \gf_{f\circ F} |_{\Sigma} \right) .
$$
Proposition~\ref{p:atlevel} can now be easily obtained from the variational principle for entropy.
For convenience, we appeal to the following result of Bowen, 
\cite[Corollary~18]{Bowen.71}, 
applied with $X=D(\Sigma)$, $Y = [0,1]$, $\phi_t = \gf_K^t|_{D(\Sigma)}$ 
and $\pi \colon D(\Sigma) \to [0,1]$, $\pi(x)=s$ if $x\in s \1 \Sigma$.

\m \ni
{\it
Let $X$, $Y$ be compact metric spaces and $\phi_t \colon X \to X$ a continuous flow.
Suppose that $\pi \colon X \to Y$ is a continuous map such that $\pi \circ \phi_t = \pi$.
Then
$$
h_{\top}(\gf) \,=\,
\sup_{y \in Y} h_{\top} \left( \phi |_{\pi^{-1}(y)} \right) .
$$ 
}
\proofend

\section{Proof of Corollary~1} \label{s:c1}

\ni
%
%
%
%
%
It suffices to prove Corollary~1 for the Hamiltonian function 
$K \colon T^*M \to \RR$ defined in~\eqref{def:K}.
Write $q_0=q$, and fix $q_1 \in M$ 
and $n \in \NN$.
Recall from the proof of Proposition~\ref{p:atlevel} that for $s \in [0,1]$
the flow $\gf_K^t |_{s\Sigma}$ is conjugate by $\psi_s$ to the flow $\gf_K^{\gs(s)\1 t}|_\Sigma$,
where $\gs (s) = f'(s^2) \1 s$.
We can choose the function~$f$ such that, in addition to the properties~\eqref{d:f},
we have
$$
\gs'(s) \,=\, 2s f''(s^2) \1 s + f'(s^2) \,>0\, \quad \mbox{for all }\, s > \eps .
$$
Then the orbits of $\gf_K^t|_\Sigma$ starting from $\Sigma_{q_0}$ at $t=0$
and arriving on $\Sigma_{q_1}$ at $t \le n$
are in bijection to points in $\gf_K^n \bigl( D_{q_0}(\Sigma) \bigr) \cap D_{q_1}(\Sigma)$.

Assume now that 
$$
q_1 \in V(q_0) := \bigcap_{n \ge 1} V_n(q_0) .
$$ 
By Corollary~\ref{c:full} the set $V(q_0)$ has full measure in~$M$.
Let $h$ be the constant from Theorem~\ref{t:floer},
that depends only on $(M,g)$.
For $q_1 \in V(q_0)$ we have, by~\eqref{e:nC}, that
$$
\# \bigl( \gf_K^n (D_{q_0}(\Sigma)) \cap D_{q_1}(\Sigma) \bigr) \,\ge\, e^{h \1 n} 
$$
provided that $n \ge N$ is large enough.
The first assertion of Corollary~1 follows.

\s
Choose a sequence $\{q_{1,j}\} \subset V(q_0)$ such that $\lim_{j \to \infty} q_{1,j} =q_1$.
Let $N \in \NN$ as chosen in the beginning of~\ref{s:Floer.growth}.
By~\eqref{e:nC}, for each $j \in \NN$ the intersection 
$\bigl( \gf_K^N (D_{q_0}(\Sigma)) \cap D_{q_{1,j}}(\Sigma) \bigr)$
is non-empty.
We thus find a sequence of $\gf_K^t$-orbits $\gg_j$ on $\Sigma$
such that $\gg_j(0) \in \Sigma_{q_0}$ and $\gg_j(t_j) \in \Sigma_{q_{1,j}}$ for some $t_j \le N$.
By the Arzel\`a--Ascoli theorem we obtain a $\gf_K^t$-orbit $\gg$ on $\Sigma$ from  
$\Sigma_{q_0}$ to $\Sigma_{q_1}$ arriving at $t \le N$.
\proofend

\begin{remark}
{\rm
Using Floer homologies $\HF_*^{[b_n,a_n]}$ with suitable action windows,
the assertions of Corollary~1 can be improved:
{\it 
If $\pi_1(M)$ has exponential growth, 
then the exponential lower bound for $\nu_n(q,q',H)$ 
holds for all pairs $q,q'$. 
Moreover, for all $q,q' \in M$ 
there are infinitely many flow lines from $S_qM$ to $S_{q'}M$,}
see~\cite{Wullschleger}.
}
\end{remark}


\section{An example}  \label{s:example}

\ni
The following example was pointed out to us by Gabriel Paternain.
Following~\cite{Mac-Pat}, we consider the group $G={\bf Sol}$,
which is the semidirect product of $\RR^2$ with $\RR$, with
coordinates $q = (x,y,z)$ and multiplication
\begin{equation*} \label{eq:sm}
(x,y,z) \star(x',y',z') \,=\, (x+e^{z}x', y+e^{-z}y',z+z').
\end{equation*}
The map $(x,y,z) \mapsto z$ is the epimorphism ${\bf Sol} \to \R$
whose kernel is the normal subgroup~$\RR^2$.
The group ${\bf Sol}$ is isomorphic to the matrix group
$$
\left(\begin{array}{ccc}
e^z&0&x\\
0&e^{-z}&y\\
0&0&1\\
\end{array}\right).
$$
It admits cocompact lattices.
Indeed, let $A \in \SL(2,\mathbb Z)$ be such that there is $P \in \GL(2,\mathbb R)$ with
$$
PAP^{-1} \,=\, \left(\begin{array}{cc}
\lambda &0\\
0&1/\lambda\\
\end{array}\right)
$$
and $\lambda>1$.
The injective homomorphism
$$
\mathbb Z^2\ltimes_{A}\mathbb{Z}\hookrightarrow \mathbf{Sol}
$$
given by $(m,n,l) \mapsto (P(m,n), (\log\lambda)\,l)$ defines
a cocompact lattice $\Gamma$ in~$\mathbf{Sol}$.
The closed 3-manifold $M := \Gamma \setminus {\bf Sol}$ is a 2-torus
bundle over the circle with hyperbolic gluing map~$A$.
Since the group~$\Gamma$ has exponential growth, $M$~is energy hyperbolic.

If we denote by $p_x$, $p_y$ and $p_z$ the momenta that
are canonically conjugate to $x$, $y$ and $z$, respectively,
then the functions
\begin{equation*} \label{eq:mom}
\begin{array}{lcl}
M_x &=& e^{z} p_{x}, \\
M_y &=& e^{-z} p_{y}, \\
M_z &=& p_{z}
\end{array}
\end{equation*}
are left-invariant functions on $T^*{\bf Sol}$. The 1-form~$\theta$ defined by
$\theta_q = e^{-z}dx$ is also left-invariant. 
The Hamiltonian on ${\bf Sol}$ defined by
\begin{equation} \label{eq:H}
2H \,=\, e^{2z}(p_{x}+e^{-z})^2+e^{-2z}p_{y}^2+p_{z}^2 \,=\, (M_x+1)^2+M_y^2+M_z^2
\end{equation}
is left-invariant and hence descends to~$M$.
Observe that $H(q,p) = \frac{1}{2} \left| p+\theta_q \right|^2$,
where the norm is induced by the left-invariant Riemannian metric 
$$
ds^2 \,=\, e^{-2z}dx^2+e^{2z}dy^2+dz^2
$$
on ${\bf Sol}$.
Note that $\left| \theta_q \right|=1$ at each $q \in M$,
and that for each $k>0$ the energy level $\Sigma_k := H^{-1}(k)$
is a $2$-sphere bundle over the graph of~$\theta$.
In terms of the functions $M_x, M_y, M_z$ we have
$\Sigma_k = M \times S_k$, where $S_k = \left\{ (M_x+1)^2+M_y^2+M_z^2 \,=\, 2k \right\}$.
Since $S_k$ encloses the origin iff $k > 1/2$,
the hypersurface~$\Sigma_k$ is fiberwise starshaped with respect to the origin 
iff $k > 1/2$.
Denote by~$\gf_H$ the Hamiltonian flow of~$H$. 
By Theorem~1, $h_{\top} (\gf_H |_{\Sigma_k}) >0$ if $k > \frac 12$.
The following proposition shows that the assumption in Theorem~1 that $\Sigma$ is fiberwise starshaped 
with respect to the origin can, in general, not be omitted.

\begin{proposition}
$h_{\top}(\gf_H |_{\Sigma_k}) >0$ if and only if $k>1/2$.
\end{proposition}

\proof
The Hamiltonian vector field of $H$ is given by
\begin{equation} \label{eq:XH}
X_H = \left\{ \begin{array}{lclclcl}
\dot{x}   &=& (M_{x}+1)e^z,          & \hspace{10mm} & \dot{M_{x}}   &=& M_{x}M_{z},\\
\dot{y} &=& M_{y}\,e^{-z}, &            & \dot{M_{y}} &=& -M_{y}M_{z},\\
\dot{z} &=& M_{z}, &           & \dot{M_{z}} &=& M_{y}^2 - M_{x}(M_{x}+1).
\end{array}
\right.
\end{equation}
The equations on the right hand side describe the Euler vector field associated to~$X_H$. 

Fix $k>0$, and abbreviate $\gf = \gf_H |_{\Sigma_k}$.
By the variational principle for entropy,
$$ 
h_{\top}(\gf) \,=\, \sup_{\mu \in \mathcal M(\gf)} h_\mu(\gf) 
$$
where $\mathcal M(\gf)$ is the set of $\gf$-invariant Borel probability measures on~$\Sigma_k$, 
and $h_\mu(\gf)$ is the entropy of~$\gf$ with respect to the measure~$\mu$. 
By the Margulis--Ruelle inequality, in turn, 
\begin{equation} \label{Ruelle}
h_\mu(\gf) \,\le\, \int_{\Sigma_k} \chi_+(q,p)\,d\mu(q,p),
\end{equation}
where $\chi_+(q,p)$ is the sum of the positive Lyapunov exponents of $(q,p) \in \Sigma_k$. 
As explained in~\cite[Section~3.2]{BP}, the Lyapunov exponents of~$\gf$ 
can be computed using the projection ${\bf Sol} \overset{z}{\longrightarrow} \RR$, 
and it turns out that
\begin{equation} \label{eq:LE}
\chi_+(q,p) \,=\, \lim_{T \to \infty} \frac{1}{T} \, \biggl| \int_0^T M_z(\gf^t(q,p))\,dt \,\biggr| .
\end{equation}
We claim that for $k \le 1/2$ the right hand side of~\eqref{Ruelle} vanishes 
for every invariant measure~$\mu$. 
Indeed, if $k<1/2$ equation~\eqref{eq:H} implies that $M_x < 0$. 
This and~\eqref{eq:XH} imply that $M_z = \dot{M_x}/M_x = \frac{d}{dt}\ln(-M_x)$, 
which, by~\eqref{eq:LE}, implies that $\chi_+(q,p)=0$ for every $(q,p) \in \Sigma_k$.

Consider now the case $k=1/2$. 
By~\eqref{eq:H} we have that $M_x(q,p) \leq 0$ for every $(q,p) \in \Sigma_k$. 
By the previous argument, if $\sup_{t \in \RR} M_x(\gf^t(q,p)) \neq 0$ then $\chi_+(q,p)=0$. 
Suppose that $(q,p)$ satisfies $\sup_{t \in \RR} M_x(\gf^t(q,p))=0$. 
It follows from~\eqref{eq:XH} that $M_xM_y$ is a first integral of $H$.
Since $\sup_{t \in \RR} M_x(\gf^t(q,p))=0$, this integral must vanish. 
Hence $M_x(q,p)=0$ or $M_y(q,p)=0$. 
If $M_x(q,p)=0$, equation~\eqref{eq:H} yields $M_y(q,p)=0$. 
By~\eqref{eq:XH}, $M_y(\gf^t(q,p))=0$ for all~$t$. 
Inspection of~\eqref{eq:XH} shows that an orbit of the Euler vector field in the circle given by 
the intersection of the plane $M_y=0$ with the sphere $(M_x+1)^2 + M_y^2 + M_z^2 = 1$ 
is either the singularity $M_x=M_y=M_z=0$ or a regular orbit that converges to the origin. 
But, by Poincar\'e Recurrence Theorem, $\mu$-almost every point $(q,p)$ is recurrent. 
Consequently only the points with $M_x=M_y=M_z=0$ may contribute to the integral in~\eqref{Ruelle}. 
These points, by~\eqref{eq:LE}, do not have positive Lyapunov exponents. 
We conclude that $h_{\top}(\gf) =0$ if $k \le 1/2$.

Recall that Theorem~1 implies $h_{\top}(\gf) >0$ for $k > 1/2$.
We give a different proof for this example:
For $k>1/2$, the Euler vector field has the singularities 
$p_\pm$ given by $M_x=M_y=0$ and $M_z=\pm\sqrt{2k-1}$. 
The measure $\mu_M = dx \wedge dy \wedge dz$ is smooth on $M \times \{p_+\}$ and $\gf$-invariant.
Moreover, it is easy to see from~\cite[Section~3.2]{BP} that the positive Lyapunov exponents of points in $M \times \{p_+\}$
are given by the magnetic flow restricted to $M \times \{p_+\}$.
Let $\mu$ be the induced measure on $\Sigma_k$ given by $\mu(A) = \mu_M(A \cap (M \times \{p_+\}))$ for every Borel set $A$.
It follows from the definition of entropy that $h_{\mu}(\gf) = h_{\mu_M}(\gf|_{M \times \{p_+\}})$.
By the Variational Principle, Pesin's formula, and \eqref{eq:LE},
$$ 
h_{\top}(\gf) \,\ge\, h_{\mu}(\gf) = h_{\mu_M}(\gf|_{M \times \{p_+\}}) \,=\, \int_{\Sigma_k} \chi_+(q,p)\,d\mu(q,p) \,=\, \sqrt{2k-1} \,>\, 0, 
$$
finishing the proof of the proposition.
\proofend


\begin{thebibliography}{99}

\bibitem{AM}
A.~Abbondandolo and P.~Majer. Lectures on the Morse complex for
infinite-dimensional manifolds. Morse theoretic methods in nonlinear
analysis and in symplectic topology, 1--74, {\it NATO Sci.\ Ser.\ II
Math.\ Phys.\ Chem.},~{\bf 217}, Springer, Dordrecht, 2006.


\bibitem{AS}
A.\ Abbondandolo and M.\ Schwarz.
On the Floer homology of cotangent bundles.
{\it Comm. Pure Appl. Math.}~{\bf 59} (2006) 254--316.


\bibitem{Benci}
V.\ Benci. 
Periodic solutions of Lagrangian systems on a compact manifold.
{\it J.~Diff.~Eq.} {\bf 63} (1986) 135--161.


\bibitem{BPS}
P.\ Biran, L.\ Polterovich and D.\ Salamon.
Propagation in Hamiltonian dynamics and relative symplectic
homology.
{\it Duke Math.\ J.} {\bf 119} (2003) 65--118.


\bibitem{Bowen.71}
R.~Bowen.
Entropy for group endomorphisms and homogeneous spaces.  
{\it Trans. Amer. Math.~Soc.}~{\bf 153} (1971) 401--414.


\bibitem{BP} 
L.~Butler and G.~Paternain. 
Magnetic flows on {\bf Sol}-manifolds: dynamical and symplectic aspects.
{\it Comm.~Math.\ Phys.}~{\bf 284} (2008) 187--202.

 

\bibitem{CFH}
K.\ Cieliebak, A.\ Floer and H.\ Hofer.
Symplectic homology. II. A general construction.
{\it Math.~Z.}~{\bf 218} (1995) 103--122.

\bibitem{CoH}
V.~Colin and Ko Honda. 
Reeb vector fields and open book decompositions.
arXiv:0809.5088 


\bibitem{Dina}
E.~I.~Dinaburg.
A connection between various entropy characterizations of
dynamical systems.
{\it Izv.~Akad. Nauk SSSR Ser.~Mat.} {\bf 35} (1971) 324--366.


\bibitem{EKP}
Ya.~Eliashberg, Sang Seon Kim, and L.~Polterovich.
Geometry of contact transformations and domains: orderability versus squeezing.
{\it Geom. Topol.}~{\bf 10} (2006) 1635--1747.


\bibitem{F2}
A.\ Floer.
A relative Morse index for the symplectic action.
{\it Comm.\ Pure Appl.\ Math.} {\bf 41} (1988) 393--407.

\bibitem{F3}
A.\ Floer.
The unregularized gradient flow of the symplectic action.
{\it Comm.\ Pure Appl.\ Math.} {\bf 41} (1988) 775--813.

\bibitem{F1}
A.\ Floer.
Morse theory for Lagrangian intersections.
{\it J.~Differential Geom.} {\bf 28} (1988) 513--547.

\bibitem{F4}
A.\ Floer.
Witten's complex and infinite-dimensional Morse theory.
{\it J.~Differential Geom.} {\bf 30} (1989) 207--221.


\bibitem{F5}
A.\ Floer. Symplectic fixed points and holomorphic spheres. {\it
Comm.\ Math.\ Phys.}~{\bf 120} (1989) 575--611.


\bibitem{FH}
A.~Floer and H.~Hofer. Symplectic homology.~I. Open sets in $C^n$.
{\it Math.~Z.}~{\bf 215} (1994) 37--88.


\bibitem{FHS}
A.~Floer, H.~Hofer and D.~Salamon.
Transversality in elliptic Morse theory for the symplectic action.
{\it Duke Math.~J.}~{\bf 80} (1995) 251--292.


\bibitem{FS:GAFA}
U.\ Frauenfelder and F.\ Schlenk.
Volume growth in the component of the Dehn--Seidel twist.
{\it Geom. Funct. Anal.}~{\bf 15} (2005) 809--838.



\bibitem{FS3}
U.\ Frauenfelder and F.\ Schlenk.
Fiberwise volume growth via Lagrangian intersections.
{\it J.~Symplectic Geom.}~{\bf 4} (2006) 117--148.


\bibitem{FS}
U.\ Frauenfelder and F.\ Schlenk. 
Hamiltonian dynamics on convex symplectic manifolds. {\it Israel J.~Math.}~{\bf 159} (2007) 1--56.


\bibitem{Gromov78}
M.\ Gromov.
Homotopical effects of dilatation.
{\it J.~Differential Geom.}~{\bf 13} (1978) 303--310.



\bibitem{Heistercamp} 
M.~Heistercamp.
The Weinstein conjecture with multiplicities on spherizations.
PhD Thesis Universit\'e de Neuch\^atel.
In preparation.


\bibitem{HMS} 
M.~Heistercamp, L.~Macarini and F.~Schlenk.
Energy surfaces in $\RR^{2n}$ and in cotangent bundles --
convex versus starshaped.
In preparation.




\bibitem{Katok.IHES}
A.~Katok.
Lyapunov exponents, entropy and periodic orbits for diffeomorphisms.  
{\it Inst. Hautes Études Sci. Publ. Math.}~{\it 51} (1980) 137--173. 


\bibitem{Katok.ETDS}
A.~Katok.
Entropy and closed geodesics.  
{\it Ergodic Theory Dynam. Systems}~{\bf 2} (1982) 339--365 (1983)


\bibitem{KH}
A.\ Katok and B.\ Hasselblatt.
{\it Introduction to the modern theory of dynamical systems.}
Encyclopedia of Mathematics and its Applications {\bf 54}.
Cambridge University Press, Cambridge, 1995.


\bibitem{KS}
M.\ Khovanov and P.\ Seidel.
Quivers, Floer cohomology, and braid group actions.
{\it J.~Amer.\ Math.\ Soc.} {\bf 15} (2002) 203--271.


\bibitem{Kli}
W.~Klingenberg.
Riemannian geometry. Second edition.
{\it de~Gruyter Studies in Mathematics,~1}.
de~Gruyter, Berlin, 1995.


\bibitem{MMP}
L.~Macarini, W.~Merry and G.\ Paternain.
In preparation.



\bibitem{Mac-Pat}
L.~Macarini and G.~Paternain.
On the stability of Ma\~n\'e critical hypersurfaces. 
to appear in {\it Calc.\ Var.\ Partial Differential Equations.} 


\bibitem{Mane}
R.~Ma\~n\'e.
On the topological entropy of geodesic flows.  
{\it J.~Differential Geom.}~{\bf 45} (1997) 74--93.


\bibitem{Mi}
J.\ Milnor.
A note on curvature and fundamental group.
{\it J.\ Differential Geometry} {\bf 2} (1968) 1--7.


\bibitem{Newh}
S.\ Newhouse.
Entropy and volume.
{\it Ergodic Theory Dynam.\ Systems} {\bf 8$\sp *$} (1988),
Charles Conley Memorial Issue, 283--299.


\bibitem{Niche}
C.~Niche.
Topological entropy of a magnetic flow and the growth of the number of trajectories.
{\it Discrete Contin.~Dyn.~Syst.}~{\bf 11} (2004) 577--580.


\bibitem{Pa1}
G.\ Paternain.
Topological entropy for geodesic flows on fibre bundles over
rationally hyperbolic manifolds.
{\it Proc.\ Amer.\ Math.\ Soc.} {\bf 125} (1997) 2759--2765.


\bibitem{Pa}
G.\ Paternain.
{\it Geodesic flows}. 
Progress in Mathematics {\bf 180}. Birkh\"auser Boston,
Inc., Boston, MA, 1999.


\bibitem{PaPa}
G.~Paternain and M.~Paternain.
Topological entropy versus geodesic entropy. 
{\it Internat. J.~Math.}~{\bf 5} (1994) 213--218. 




\bibitem{PP2}
G.\ Paternain and J.\ Petean.
Zero entropy and bounded topology.
{\it Comment. Math.~Helv.}~{\bf 81} (2006) 287--304.


\bibitem{Sa-90}
D.\ Salamon.
Morse theory, the Conley index and Floer homology.
{\it Bull.\ London Math.\ Soc.}~{\bf 22} (1990) 113--140.


\bibitem{SZ}
D.~Salamon and E.~Zehnder. Morse theory for periodic solutions of
Hamiltonian systems and the Maslov index. {\it Comm.\ Pure Appl.\
Math.}~{\bf 45} (1992) 1303--1360.


\bibitem{Schw} 
M.~Schwarz.
{\it Morse homology}. Progress in Mathematics~{\bf 111}.
Birkh\"auser Verlag, Basel, 1993.

\bibitem{Seidel.biased}
P.~Seidel. 
A biased view of symplectic cohomology. 
{\it Current Developments in Mathematics}~{\bf 2006} (2008) 211--253. 


\bibitem{V2} 
C.\ Viterbo.
Functors and computations in Floer homology with applications. I.
{\it Geom.\ Funct.\ Anal.} {\bf 9} (1999) 985--1033.


\bibitem{Wullschleger} 
R.\ Wullschleger.
Slow entropy of Reeb flows in spherizations.
PhD Thesis Universit\'e de Neuch\^atel.
In preparation.


\bibitem{Y} 
Y.\ Yomdin.
Volume growth and entropy.
{\it Israel J.~Math.} {\bf 57} (1987) 285--300.

\end{thebibliography}
\enddocument